\documentclass[twoside, 11pt, final]{article} 

\usepackage[utf8]{inputenc}
\usepackage[T1]{fontenc}
\usepackage{amsmath, amssymb, amsthm, amstext, amsfonts}
\usepackage{dsfont}
\usepackage[svgnames]{xcolor}
\usepackage{graphics,graphicx,subfigure}
\usepackage{empheq}

\usepackage{a4wide}
\usepackage{hyperref}
\usepackage{enumerate}

\usepackage{pgf,tikz}

\numberwithin{equation}{section}
\numberwithin{figure}{section}
\allowdisplaybreaks[4]

\newtheorem{proposition}{Proposition}[section]
\newtheorem{theorem}[proposition]{Theorem}
\newtheorem{lemma}[proposition]{Lemma}
\newtheorem{definition}[proposition]{Definition}
\newtheorem{corollary}[proposition]{Corollary}

\newenvironment{proofof}[1]{\medskip\noindent{\textbf{Proof~of~#1.}}%
  \hspace{1pt}}{\hspace{-5pt}{\nobreak\quad\nobreak\hfill\nobreak%
    $\square$\vspace{2pt}\par}\smallskip\goodbreak}

\newcommand{\fonction}[5]{
    \begin{array}{l@{\,}ccc}
      #1\colon
      & #2
      & \longrightarrow
      & #3
      \\
      & \displaystyle{#4}
      & \longmapsto
      & \displaystyle{#5}
    \end{array}}

\newcommand{\limit}[2]{{\ \underset{#1 \to #2}{\longrightarrow} \ }}

\renewcommand{\d}[1]{\mathinner{\mathrm{d}{#1}}}
\renewcommand{\epsilon}{\varepsilon}
\newcommand{\sgn}{\mathop{\rm sgn}}

\newcommand{\N}{\mathbb{N}}

\newcommand{\R}{\mathbb{R}}

\newcommand{\C}[1]{\mathbf{C}^{#1}} \newcommand{\Cc}[1]{\mathbf{C}_\mathbf{c}^{#1}}
\renewcommand{\L}[1]{\mathbf{L}^{#1}} \newcommand{\Lloc}[1]{\mathbf{L}_{\mathbf{loc}}^{#1}}
\newcommand{\W}[2]{\mathbf{W}^{#1, #2}}
\newcommand{\Lip}{\mathbf{Lip}} 
\newcommand{\BV}{\mathbf{BV}}
\newcommand{\BVloc}{\mathbf{BV}_{\mathbf{loc}}}

\newcommand{\modulo}[1]{{\left|#1\right|}}
\newcommand{\norma}[1]{{\left\|#1\right\|}}

\newcommand{\J}{\mathcal{J}}

\newcommand{\cG}{\mathcal{G}}

\newcommand{\cR}{\mathcal{R}}

\newcommand{\bK}{\mathbf{K}}

\newcommand{\flow}{\mathcal{F}}
\renewcommand{\u}{z}

\begin{document}

\title{\textbf{Initial Data Identification in Space Dependent \\
    Conservation Laws and Hamilton-Jacobi Equations}}

\author{Rinaldo M. Colombo$^1$ \qquad Vincent Perrollaz$^2$ \qquad
  Abraham Sylla$^3$}

\date{\today}

\maketitle

\footnotetext[1]{INdAM Unit \& Department of Information Engineering,
  University of Brescia, Italy.}

\footnotetext[2]{Institute Denis Poisson, University of Tours, CNRS
  UMR 7013, University of Orl\'eans, France.}

\footnotetext[3]{Department of Mathematics and Applications,
  University of Milano -- Bicocca, Italy.}
\maketitle

\begin{abstract}
  \noindent Consider a Conservation Law and a Hamilton-Jacobi equation
  with a flux/Hamiltonian depending also on the space variable. We
  characterize first the attainable set of the two equations and,
  second, the set of initial data evolving at a prescribed time into a
  prescribed profile. An explicit example then shows the deep
  differences between the cases of $x$-independent and $x$-dependent
  fluxes/Hamiltonians.

  \smallskip


\end{abstract}


\section{Introduction}

We characterize the \emph{inverse designs} for Conservation Laws and
for Hamilton-Jacobi equations. They are the sets of those initial data
that, separately for the two equations, evolve into a given profile
after a given positive time.

As is well known, both Conservation Laws and Hamilton-Jacobi equations
generate Lipschitz continuous semigroups whose orbits are solutions,
either in the entropy sense or in the viscosity sense. However, the
insurgence of singularities implies that these evolutions may not be
time reversible, in general. As a result, inverse designs, when non
empty, may well display interesting --- infinite dimensional ---
geometric or topological properties.

From a control theoretic point of view, the characterization of
inverse designs solves the most elementary controllability problem,
thus playing a key role in subsequent developments. Indeed, the first
step in the study of inverse designs consists in a full
characterization of the attainable sets, i.e., of the profiles leading
to non empty inverse designs. In this connection, the current
literature offers a few results, typically limited to the
$x$-independent case.  We refer the reader to~\cite{MR1616586} for a
characterization of the attainable set for a conservation law (here,
with boundary); to~\cite{MR4503821} for a result on the attainable set
for Hamilton-Jacobi equations in several space dimensions and
to~\cite{MR3489384} for the case of an $x$-dependent source term. A
triangular system of conservation laws is considered
in~\cite{MR3394696}.

\smallskip

Below, we proceed beyond reachable sets and fully characterize inverse
designs.

\smallskip

\noindent More precisely, we consider the conservation law
\begin{equation}
  \label{eq:cl}\tag{\textbf{CL}}
  \left\{
    \begin{array}{r@{\,}c@{\,}l@{\qquad}r@{\,}c@{\,}l}
      \partial_t u + \partial_x \left(H(x, u)\right)
      & =
      & 0
      & (t, x)
      & \in
      & \mathopen]0, +\infty\mathclose[ \times \R
      \\
      u(0, x)
      & =
      & u_o(x)
      & x
      & \in
      & \R
    \end{array}
  \right.
\end{equation}
and the Hamilton-Jacobi equation
\begin{equation}
  \label{eq:hj}\tag{\textbf{HJ}}
  \left\{
    \begin{array}{r@{\,}c@{\,}l@{\qquad}r@{\,}c@{\,}l}
      \partial_t U + H(x, \partial_x U)
      & =
      & 0
      & (t, x)
      & \in
      & \mathopen]0, +\infty \mathclose[ \times \R \nonumber
      \\
      U(0, x)
      & =
      & U_o(x)
      & x
      & \in
      & \R
    \end{array}
  \right.
\end{equation}
both in the scalar, one dimensional, \emph{non homogeneous}, i.e.,
$x$-dependent, case. Denote by
\begin{equation}
  \label{eq:1}
  S^{CL} \colon \R_+ \times \L{\infty}(\R; \R) \to \L{\infty}(\R; \R)
  \quad \mbox{ and } \quad
  S^{HJ} \colon \R_+ \times \Lip(\R; \R) \to \Lip(\R; \R) \,,
\end{equation}
respectively, the semigroups whose orbits are entropy solutions
to~\eqref{eq:cl} and viscosity solutions to~\eqref{eq:hj},
see~\cite[\S~2.5]{CPS2022}. For any positive $T$ and for any assigned
profiles $w \in \L{\infty}(\R; \R)$ and $W \in \Lip(\R; \R)$, the
inverse designs are
\begin{equation}
  \label{eq:2}
  \begin{array}{rcl}
    I_T^{CL}(w)
    & \coloneqq
    & \left\{u_o \in \L{\infty}(\R; \R) \colon S_T^{CL}u_o = w \right\}
      \quad \mbox{ and}
    \\
    I_T^{HJ}(W)
    & \coloneqq
    & \left\{U_o \in \Lip(\R; \R) \colon S_T^{HJ} U_o = W \right\}.
  \end{array}
\end{equation}

In the homogeneous --- $x$-independent --- case, a general
characterization of $I^{CL}_T (w)$ and $I^{HJ}_T (W)$ is given
in~\cite{CP2020}. Other more specific results in this setting
are~\cite{LZ2021}, devoted to Burgers' equation; \cite{MR1616586},
specific to boundary value problems arising in the modeling of
vehicular traffic. The multi--dimensional setting is considered
in~\cite{EZ2020}, specifically in the case of~\eqref{eq:hj}. A
classical reference for analytic techniques used in these papers
is~\cite{MR1722801}.

The present non homogeneous case significantly differs from the
homogeneous one and significantly less results in the literature are
available. The explicit example constructed below shows that when $H$
depends on $x$ (even smoothly), the inverse design $I^{CL}_T (w)$ may
have properties in a sense opposite to the general ones that hold in
the homogeneous case, according to~\cite{CP2020}. In particular, for
instance, the results in~\cite{CP2020} ensure that in the
$x$-independent case
\begin{displaymath}
  S_T^{CL} \left(\L\infty (\R;\R)\right)
  =
  {\mathop\mathbf{cl}}_{\strut\L1}
  \left(
    S_T^{CL}
    \left\{u_o \in \C1 (\R;\R) \colon S_tu_o \in \C1 (\R;\R)
      \mbox{ for all } t \in [0,T]\right\}
  \right)
\end{displaymath}
which can be false when $H$ depends on $x$, as in the case of the
example in Section~\ref{sec:pec}. It thus appears that non homogeneous
Conservation Laws are, in a sense, \emph{more singular} than
homogeneous ones.

Assume $I_T^{CL} (w)$ is non empty. Then, in the $x$-independent case,
the presence of a shock in $w$ is a necessary and sufficient condition
for $I^{CL}_T (w)$ to be infinite or, equivalently, $I^{CL}_T (w)$ is
a singleton if and only if $w$ is continuous. More precisely, in the
$x$-independent case, the presence of a shock in $w$ implies that
$I^{CL}_T (w)$ is a close convex cone without extremal faces of finite
dimension. On the contrary, in the $x$-dependent case, we exhibit an
example where $I^{CL}_T (w)$ is a singleton although $w$ displays a
shock. This is explained in Section~\ref{sec:pec}, where the theory of
generalized characteristics, see~\cite{Dafermos1977}, is deeply
exploited.
\begin{figure}[!ht]
  \begin{center}
    \includegraphics[scale=0.5]{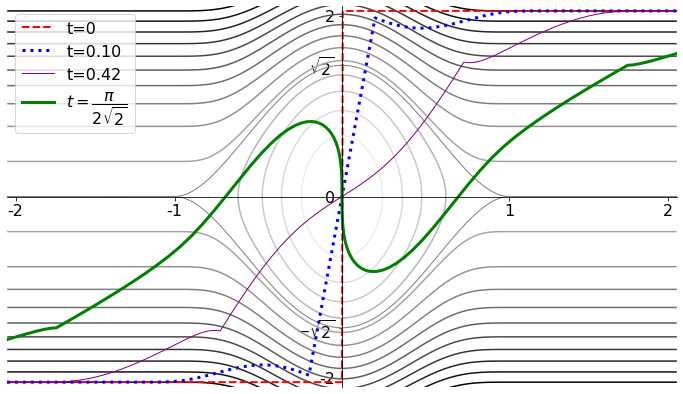}
  \end{center}
  \vspace*{-5mm}
  \caption{Superposition of a solution to~\eqref{eq:cl} at different
    times with the orbits of the Hamiltonian
    system~\eqref{eq:ode_system}. $x$ (or $q$) is on the horizontal
    axis and $u$ (or $p$) on the vertical axis. As proved later in
    Theorem~\ref{th:counter_example}, the initial datum~\eqref{eq:12}
    is the unique one that evolves into the depicted profiles where,
    at time $T = \left.\pi \middle/ (2\sqrt2)\right.$, a shock
    arises.}
  \label{fig:nice}
\end{figure}
Graphs of the constructed solution are in Figure~\ref{fig:nice}.

We defer further remarks on these differences to
Theorem~\ref{th:counter_example} and to the subsequent discussion. Let
us recall that a first step in this direction, limited to the study of
the attainable set, is~\cite{MR4188826}, where $H$ in~\eqref{eq:cl}
consists of an expression for $x>0$ and another expression for $x<0$,
see also the related preprint~\cite{AdimurthiGhoshal2020}.

The analytic techniques developed below take advantage of the deep
connection between~\eqref{eq:cl} and~\eqref{eq:hj}. We know, on the
basis of~\cite{CPS2022}, that both these Cauchy problems are
(globally) well posed under the same set of assumptions, namely
\begin{align}
  \label{eq:smoothness}\tag{\textbf{C3}}
  \mbox{\textbf{Smoothness:}}
  & \;\, H \in \C{3}(\R^2;\R) \,.
  \\[6pt]
  \label{eq:space_dpcy}\tag{\textbf{CNH}}
  \mbox{\textbf{Compact~NonHomogeneity:}}
  & \begin{array}{l}
      \exists X > 0, \; \forall (x, p) \in \R^2,
      \\
      \modulo{x} \geq X \implies \partial_x H(x, p) = 0.
    \end{array}
  \\[6pt]
  \label{eq:convexity}\tag{\textbf{CVX}}
  \mbox{\textbf{Strong~Convexity:}}
  &\begin{array}{l}
     \forall x \in \R, \;
     p \mapsto \partial_p H(x, p) \mbox{ is an}
     \\
     \text{increasing $\C{1}-$diffeomorphism}
     \\
     \text{of $\R$ onto itself.}
   \end{array}
\end{align}
Rather than tackling directly the characterization of the inverse
design for~\eqref{eq:cl}, we do it for~\eqref{eq:hj} and use the
correspondence to get back to~\eqref{eq:cl}.

Assumption~\eqref{eq:convexity} implies that $H$ is strictly convex
with respect to the second variable. As is well known, the mappings
$x \mapsto -x$ and $H \mapsto -H$ transform the convex case into the
concave one, and \emph{vice versa}. Recall that~\eqref{eq:convexity}
is a recurrent assumption in the context of~\eqref{eq:hj} where it
allows a connection to optimal control, see~\cite{BCDBook, BarlesBook,
  CSBook}. On the contrary, the use of
Assumption~\eqref{eq:space_dpcy} in conservation laws, to the authors'
knowledge, was recently introduced in~\cite{CPS2022}.

It is worth noting that the
assumptions~\eqref{eq:smoothness}--\eqref{eq:space_dpcy}--\eqref{eq:convexity}
comprise fluxes (Hamiltonians) that do not fit in the classical {Kru\v
  zkov} paper~\cite{Kruzhkov1970}. Indeed,
following~\cite[Example~1.1]{CPS2022} consider the Hamiltonian
\begin{equation}
  \label{eq:traffic_flow1}
  H(x, u) \coloneqq V(x) \, u \, \left(1 - \frac{u}{R(x)} \right),
\end{equation}
where $V, R \in \C{3}(\R; \R)$ are both strictly positive and with
compactly supported derivative.  The conservation
law~\eqref{eq:cl}--\eqref{eq:traffic_flow1} describes the time
evolution of the density $u = u (t,x)$ of a flow of vehicles along a
one-dimensional road that allows a space dependent maximal density
$R = R(x)$ and maximal speed $V = V(x)$. It is readily checked that
$H$ in~\eqref{eq:traffic_flow1} satisfies~\eqref{eq:smoothness}, and
it is strongly concave --- analogously to~\eqref{eq:convexity}. On the
other hand, this $H$ may not meet the assumptions
of~\cite{Kruzhkov1970}. In particular, it fails the growth assumption
$\sup_{(x, u) \in \R^2} \left(-\partial_{xu}^2 H(x, u) \right) <
+\infty$, see~\cite[Formula~(4.2)]{Kruzhkov1970}.

While inverse design refers to going back in time, the \emph{dual}
approach is connected to the problem of the compactness of the range
of the semigroup $S_t^{CL}$, apparently considered only in the
homogeneous case~\cite{zbMATH02183964}, extended
in~\cite{zbMATH06067762} to balance laws, but the case of fluxes
depending on the space variable is, to our knowledge, still open.

\smallskip

The next section provides the basic background. Then, on the basis
of~\cite{CPS2022}, Section~\ref{sec:ext} extends to the $x$-dependent
case several classical results, see~\cite{CP2020}. On the contrary,
the example constructed in Section~\ref{sec:pec} shows how deep can be
the differences between the homogeneous and non homogeneous case. All
proofs are deferred to Section~\ref{sec:inverse_design}.

\section{Notations and Definitions}
\label{sec:main_results}

Recall the classical definition of entropy
solution~\cite[Definition~1]{Kruzhkov1970}, as tweaked
in~\cite{CPS2022}.

\begin{definition}
  \label{def:entropy_solution}
  Fix $u_o \in \L{\infty}(\R;\R)$. A bounded function
  $u \in \L{\infty}(\R_+ \times \R; \R)$ is a \emph{solution}
  to~\eqref{eq:cl} if for all test functions
  $\phi \in \Cc{\infty}(\R_+ \times \R; \R_+)$ and for all scalar
  $k \in \R$:
  \begin{displaymath}
    \begin{array}{rcl}
      \displaystyle
      \int_0^{+\infty} \!\! \int_\R
      \modulo{u (t,x) - k} \, \partial_t \phi (t,x)
      \d{x} \d{t}
      \\
      \displaystyle
      +
      \int_0^{+\infty} \!\! \int_\R
      \sgn \left(u (t,x) - k\right) \;
      \left(H \left(x,u (t,x)\right) - H (x,k)\right)
      \, \partial_x \phi (t,x)
      \d{x} \d{t}
      \\
      \displaystyle
      -
      \int_0^{+\infty} \int_\R
      \sgn\left(u (t,x) - k\right) \, \partial_x H (x,k) \, \phi (t,x)
      \d{x} \d{t}
      \\
      \displaystyle
      + \int_\R \modulo{u_o(x) - k} \, \phi(0,x) \d{x}
      & \geq
      &  0 \,.
    \end{array}
  \end{displaymath}
\end{definition}

Definition~\ref{def:entropy_solution}, taken from
by~\cite[Definition~2.1]{CPS2022} is apparently weaker than the
classical {Kru\v zkov} definition since it does not require the
\emph{``trace at $0$ condition''}~\cite[Formula~(2.2)]
{Kruzhkov1970}. Nevertheless, under Assumption~\eqref{eq:smoothness},
Definition~\ref{def:entropy_solution} ensures uniqueness and uniform
$\Lloc1$--continuity in time of the solution, as proved
in~\cite[Theorem~2.6]{CPS2022}.

The following Lemma ensures the existence of left and right traces in
the space variable at any point. In the homogeneous ---
$x$-independent --- case, this is classically obtained through the
well known Oleinik estimates~\cite[Theorem~11.2.1 and
Theorem~11.2.2]{DafermosBook}.

\begin{lemma}
  \label{lmm:regularity_profile_cl}
  Let $H$ satisfy~\eqref{eq:smoothness}, \eqref{eq:space_dpcy}
  and~\eqref{eq:convexity}. Fix $T>0$ and $w \in \L{\infty}(\R; \R)$
  so that $I_T^{CL}(w) \neq \emptyset$. Then, for all $x \in \R$, $w$
  admits finite left and right traces at $x$.
\end{lemma}

\noindent The proof is deferred to
Section~\ref{sec:inverse_design}. Once this Lemma is proved, we are
able to use Dafermos' techniques based on generalized characteristics
from~\cite{Dafermos1977}, where solutions are however
\textbf{required} to have traces at each point. Alternatively, another
reference is~\cite[Chapter~10]{DafermosBook}
or~\cite[Section~11.11]{DafermosBook} for the inhomogeneous case, but
here solutions are \textbf{required} to be in $\BV$. Thus, particular
care has to be taken here to avoid circular arguments.

\smallskip

We now recall the framework of viscosity solutions to~\eqref{eq:hj},
introduced by Crandall--Lions.

\begin{definition}[{\cite[Definition~5.3]{CL1983}}]
  \label{def:viscosity_solution}
  Let $U \in \Lip([0,T] \times \R; \R)$ satisfy $U (0) = U_o$.
  \begin{enumerate}[\rm(i)]
  \item\label{item:def:subs} $U$ is a subsolution to~\eqref{eq:hj}
    when for all test functions
    $\phi \in \C1 (\mathopen]0,T\mathclose[ \times \R ; \R)$ and for
    all $(t_o, x_o) \in \mathopen]0,T\mathclose[ \times \R$, if
    $U - \phi$ has a point of local maximum at the point $(t_o, x_o)$,
    then
    $\partial_{t} \phi(t_o, x_o) + H\left(x_o,\partial_x \phi(t_o,
      x_o)\right) \leq 0$;

  \item\label{item:def:sups} $U$ is a supersolution to~\eqref{eq:hj}
    when for all test functions
    $\phi \in \C1 (\mathopen]0,T\mathclose[ \times \R ; \R)$ and for
    all $(t_o, x_o) \in \mathopen]0,T\mathclose[ \times \R$, if
    $U - \phi$ has a point of local minimum at the point $(t_o, x_o)$,
    then
    $\partial_{t} \phi(t_o, x_o) + H\left(x_o,\partial_x \phi(t_o,
      x_o)\right) \geq 0$.

  \item\label{item:def:visco} $U$ is a viscosity solution
    to~\eqref{eq:hj} if it is both a supersolution and a subsolution.
  \end{enumerate}
\end{definition}

\noindent The literature offers a standardized framework for the well
posedness of~\eqref{eq:cl}, typically referred to the classical
paper~\cite{Kruzhkov1970}, see also~\cite{DafermosBook}. On the
contrary, a wide variety of assumptions are available, where results
ensuring the well posedness of~\eqref{eq:hj} can be proved, see for
instance~\cite{BCDBook, BarlesBook, CSBook, CL1983} or the
textbooks~\cite[Chapter~9]{MR2347697},
\cite[Chapter~10]{EvansBook}. Here we recall in
particular~\cite{SyllaPhD}, devoted to the convex case,
and~\cite{CPS2022} where the two equations are considered under the
\textbf{same} set of assumptions, thus allowing a detailed description
of the correspondence between the solutions to the two
equations. Indeed, the orbits of the semigroups~\eqref{eq:1} are
solution to~\eqref{eq:cl} in the sense of
Definition~\ref{def:entropy_solution}, respectively~\eqref{eq:hj} in
the sense of Definition~\ref{def:viscosity_solution},
see~\cite[Theorem~2.18 and Theorem~2.19]{CPS2022}. Thanks to their
$\Lloc1$ continuity, both these semigroups a uniquely defined for all
$t \in \R_+$.

For any positive $T$ and for any assigned profiles
$w \in \L{\infty}(\R; \R)$ and $W \in \Lip(\R; \R)$, we first present
conditions ensuring that the sets $I^{CL}_T (w)$ and $I^{HJ}_T (W)$
in~\eqref{eq:2} are not empty and then prove geometrical/topological
properties. In light of the correspondence $U \to u = \partial_x U$
between $S^{HJ}$ and $S^{CL}$, see~\cite[Theorem~2.20]{CPS2022} or
also~\cite{MR2795714, MR2374224, CP2020, KarlsenRisebro2002}, each of
the two characterizations can be deduced from the other one.

As usual, in connection with~\eqref{eq:hj} and~\eqref{eq:cl}, we use
of the system of ordinary differential equations
\begin{equation}
  \label{eq:ode_system}\tag{\textbf{HS}}
  \left\{
    \begin{array}{rcl}
      \dot q
      & =
      & \partial_u H(q, p)
      \\
      \dot p
      & =
      & - \partial_x H(q, p)
    \end{array}
  \right.
\end{equation}
which we consider equipped with initial or with final
conditions. Basic properties of~\eqref{eq:ode_system}
under~\eqref{eq:smoothness}--\eqref{eq:space_dpcy}--\eqref{eq:convexity}
are proved in Lemma~\ref{lem:HS} and in the subsequent ones.  For a
fixed positive $T$, with reference to~\eqref{eq:ode_system}, we also
introduce the set
\begin{equation}
  \label{eq:hamiltonian_rays}
  \cR_T \coloneqq
  \left\{
    q \in \C{1}([0,T];\R)  \colon
    \exists \, p \in \C{1}([0,T]; \R) \text{ such that } (q, p)
    \mbox{ solves~{\eqref{eq:ode_system}}}
  \right\} \,.
\end{equation}
whose elements we call Hamiltonian rays. For all
$w \in \L{\infty}(\R; \R)$ such that $I_T^{CL}(w) \neq \emptyset$, so
that Lemma~\ref{lmm:regularity_profile_cl} applies and we can define
\begin{equation}
  \label{eq:4}
  \fonction{\pi_w}{\R}{\R}{x}{q(0),}
  \mbox{ where }
  (q,p) \mbox{ solves~\eqref{eq:ode_system}}
  \mbox{ with datum }
  \left\{
    \begin{array}{r@{\,}c@{\,}l}
      q (T)
      & =
      & x
      \\
      p (T)
      & =
      & w (x-) \,.
    \end{array}
  \right.
\end{equation}
The map $\pi_w$ assigns to $x \in \R$ the intersection of the minimal
backward characteristics emanating from $(T, x)$,
see~\cite[Definition~3.1, Theorems~3.2 and~3.3]{Dafermos1977}, with
the axis $t=0$. Lemma~\ref{lmm:regularity_profile_cl} and
Lemma~\ref{lem:HS} ensure that $\pi_w$ is well defined. Remark that in
the $x$-independent case, all Hamiltonian rays are straight lines, as
also any extremal characteristics, a key simplification exploited
in~\cite[Formula~(2.3)]{CP2020}.

As is well known, thanks to~\eqref{eq:convexity}, Hamilton-Jacobi
equation~\eqref{eq:hj} is deeply related and motivated by the search
for minima of functionals of the type
\begin{equation}
  \label{eq:8}
  \fonction{\J_{t}}{\W{1}{1}([0,t]; \R)}{\R}{y}{\int_0^{t} L\left(y(s), \dot y(s)\right) \d{s} + U_o\left(y(0)\right)}
\end{equation}
where $U_o \in \Lip(\R; \R)$ and $L$ is the Legendre transform of $H$
in $p$, i.e.,
\begin{equation}
  \label{eq:legendre}
  \fonction{L}{\R^2}{\R}{(x,v)}{\sup_{p \in \R} \left( p\, v - H(x, p) \right)\,.}
\end{equation}
As general references for this minimization problem, we refer
to~\cite[Chapter~5]{CSBook}, \cite[Part~III]{ClarkeBook},
\cite[Chapter~3]{EvansBook}. Below, for detailed proofs about the
connection between solutions to~\eqref{eq:hj} and to minimization
problems in our specific functional setting, we often refer
to~\cite[\S~8.3]{SyllaPhD}. Recall, in particular, that $U$
solves~\eqref{eq:hj} if and only if for all
$(T,x) \in \mathopen[0, +\infty\mathclose[ \times \R$,
\begin{equation}
  \label{eq:24}
  U (T,x)
  =
  \inf_{\substack{\gamma(T) = x \\ \gamma \in \cR_T}}
  \left( \int_0^T L\left(\gamma(s), \dot \gamma(s)\right) \d{s} +
    U_o\left(\gamma(0)\right) \right) \,,
\end{equation}
see~\cite[Corollary~8.3.15]{SyllaPhD}. Note moreover that
by~\cite[Theorem~8.3.12]{SyllaPhD}
\begin{displaymath}
  U (T,x)
  =
  \inf_{\substack{\gamma(T) = x \\ \gamma \in \Lip ([0,T];\R)}}
  \left( \int_0^T L\left(\gamma(s), \dot \gamma(s)\right) \d{s} +
    U_o\left(\gamma(0)\right) \right) \,,
\end{displaymath}
As a first step, we verify that the present
assumptions~\eqref{eq:smoothness}--\eqref{eq:space_dpcy}--\eqref{eq:convexity}
allow to apply the results in~\cite{CPS2022}, where convexity was
relaxed to genuine nonlinearity and uniform coercivity.

\begin{proposition}
  \label{prop:hypOk}
  Let $H$
  satisfy~\eqref{eq:smoothness}--\eqref{eq:space_dpcy}--\eqref{eq:convexity}. Then. the
  following properties hold:
  \begin{align}
    \label{eq:UC}\tag{\textbf{UC}}
    \mathbf{Uniform~Coercivity:}
    & \qquad  \begin{array}{@{}l}
                \forall \, h \in \R
                \quad
                \exists \, \mathcal{U}_h \in \R \colon
                \forall \, (x,u) \in \R^2
                \\
                \mbox{if }
                \modulo{H (x,u)} \leq h
                \mbox{ then }
                \modulo{u} \leq \mathcal{U}_h  \,.
              \end{array}
    \\[6pt]
    \label{eq:WGNL}\tag{\textbf{WGNL}}
    \mathbf{Weak~Genuine~NonLinearity:}
    & \qquad  \begin{array}{@{}l}
                \mbox{for a.e.} \, x \in \R
                \mbox{ the set }
                \\
                \left\{
                w \in \R \colon \partial^2_{ww} H (x,w) = 0
                \right\}
                \\
                \mbox{has empty interior.}
              \end{array}
  \end{align}
\end{proposition}

\noindent The proof is deferred to Section~\ref{sec:inverse_design}.

\section{Extensions from Homogeneous to Non Homogeneous}
\label{sec:ext}

This section is focused on those properties known to hold in the
homogeneous case, see~\cite{CP2020}, whose statement admits a natural
extension to the non homogeneous case. However, the proofs typically
require a new approach.

An interesting connection between~\eqref{eq:cl} and~\eqref{eq:hj} is
the following result, which shows that minimal and maximal backward
characteristics are minima of the functional~\eqref{eq:8}.

\begin{theorem}
  \label{th:two_min}
  Let $H$
  satisfy~\eqref{eq:smoothness}--\eqref{eq:space_dpcy}--\eqref{eq:convexity}. Fix
  $U_o \in \Lip(\R; \R)$ and let $U$ solve~\eqref{eq:hj} in the sense
  of Definition~\ref{def:viscosity_solution}. Fix
  $(t, x) \in \mathopen]0,+\infty\mathclose[ \times \R$ and let
  $\check\zeta$, respectively $\hat\zeta$, be the minimal,
  respectively maximal, backwards characteristics, related to
  $u = \partial_x U$ which solves~\eqref{eq:cl}, emanating from
  $(t, x)$, see~\cite[Definition~3.1]{Dafermos1977}. Then, with
  reference to the functional~\eqref{eq:8},
  \begin{displaymath}
    U(t, x) = \J_t(\check\zeta) = \J_t(\hat\zeta) \,.
  \end{displaymath}
\end{theorem}

\noindent The proof is deferred
to~\S~\ref{subsec:proof-theor-refth:tw}.

We are now ready to state the conditions ensuring that $I^{HJ}_T (W)$,
as defined in~\eqref{eq:2}, is not empty. In other words, the next
result completely characterizes the reachable set for~\eqref{eq:hj}.

\begin{theorem}
  \label{th:nonemptiness}
  Let $H$
  satisfy~\eqref{eq:smoothness}--\eqref{eq:space_dpcy}--\eqref{eq:convexity}. Fix
  $T>0$, $W \in \Lip(\R; \R)$ and define
  \begin{equation}
    \label{eq:inverse_design}
    \fonction{U_o^*}{\R}{\R}{x}{\sup_{\substack{q(0) = x \\ q \in \cR_T}}
      \left(W\left(q(T)\right) - \int_0^T L\left(q(s), \dot q(s)\right) \d{s}
      \right)\,,}
  \end{equation}
  where $L$ is as in~\eqref{eq:legendre} and $\mathcal{R}_T$ as
  in~\eqref{eq:hamiltonian_rays}. Then, the following conditions are
  equivalent:
  \begin{enumerate}[\bf(1)]

  \item\label{item:2} $U_o^* \in I_T^{HJ}(W)$.

  \item\label{item:1} $I_T^{HJ}(W) \neq \emptyset$.

  \item\label{item:3} The set
    \begin{equation}
      \label{eq:graph}
      \cG \coloneqq  \left\{ (x_o, x_T) \in \R^2 \colon
        \exists \, q \in \cR_T,
        \begin{array}{rl@{}}
          (i)
          & q(0) = x_o\,, \; q(T) = x_T
          \\
          (ii)
          & U_o^*(x_o) = W(x_T) - \int_0^T L\left(q(s), \dot q(s)\right)
            \d{s}
        \end{array}
      \right\}
    \end{equation}
    has the following property:
    \begin{equation}
      \label{eq:graph_maximal}
      (x_o, x_T'), (x_o, x_T'') \in \cG \implies
      \forall \, x_T \in [x_T', x_{T}''], \; (x_o, x_T) \in \cG \, .
    \end{equation}
  \end{enumerate}
  \noindent Moreover, any of the conditions above implies that the map
  $\pi_{W'}$ defined in~\eqref{eq:4} is well defined and
  nondecreasing.
\end{theorem}

\noindent The proof is deferred
to~\S~\ref{subsec:proofs-relat-sect}. The set $\mathcal{G}$ is more
readily interpreted from the point of view of~\eqref{eq:cl}. In
particular, \eqref{eq:graph_maximal} describes the structure of
rarefaction-like waves and, limited to the $x$-independent case,
deriving the condition on $\mathcal{G}$ from the property of
$\pi_{W'}$ is straightforward. In this connection, the $x$-dependent
case is significantly more intricate. We signal in Lemma~\ref{lmm:cs3}
additional properties of the set $\mathcal{G}$.

We are now ready to provide a full and general characterization of the
inverse designs.

\begin{theorem}
  \label{th:belonging}
  Let $H$ satisfy~\eqref{eq:smoothness}, \eqref{eq:space_dpcy}
  and~\eqref{eq:convexity}. Fix $T>0$ and $W \in \Lip(\R; \R)$ such
  that $I_T^{HJ}(W) \neq \emptyset$ and define $U_o^*$ as
  in~\eqref{eq:inverse_design}. Then, for all $U_o \in \Lip(\R; \R)$,
  \begin{equation}
    \label{eq:belonging}
    U_o \in I_T^{HJ}(W) \iff
    \left\{
      \begin{aligned}
        (i) & \; U_o \geq U_o^* \\
        (ii) & \; U_o = U_o^* \; \text{on } \overline{\pi_{W'}(\R)}.
      \end{aligned}
    \right.
  \end{equation}
  where $\pi_{W'}$ is defined as in~\eqref{eq:4}.
\end{theorem}

\noindent The proof is deferred to~\S~\ref{sec:proofs-relat-sect}.

\begin{corollary}
  \label{cor:geometrical}
  Let $H$ satisfy~\eqref{eq:smoothness}, \eqref{eq:space_dpcy}
  and~\eqref{eq:convexity}. Fix $T>0$ and $W \in \Lip(\R; \R)$ such
  that $I_T^{HJ}(W) \neq \emptyset$. Then, $I_T^{HJ}(W)$ is a closed
  convex cone with vertex $U_o^*$, defined
  in~\eqref{eq:inverse_design} and moreover
  $U_o^* = \min I_T^{HJ}(W)$.
\end{corollary}

\noindent The proof is an immediate consequence of the
characterization provided by Theorem~\ref{th:belonging}.

Corollary~\ref{cor:geometrical} admits a clear counterpart related
to~\eqref{eq:cl}, on the basis of the correspondence
between~\eqref{eq:cl} and~\eqref{eq:hj} proved
in~\cite[Theorem~2.20]{CPS2022}. An analogous characterization in the
$x$-independent case is provided by~\cite[Proposition~5.2,
Item~(G2)]{CP2020}.

\begin{corollary}
  \label{cor:geoCL}
  Let $H$ satisfy~\eqref{eq:smoothness}, \eqref{eq:space_dpcy}
  and~\eqref{eq:convexity}. Fix $T>0$ and $w \in \L\infty(\R; \R)$
  such that $I_T^{CL}(w) \neq \emptyset$. Then, $I_T^{CL}(w)$ is a
  closed convex cone with vertex $u_o^*$, defined by
  $u_o^* = \partial_x U_o^*$ and $U_o^*$ is as
  in~\eqref{eq:inverse_design}.
\end{corollary}

The latter corollary extends to the $x$-dependent case some of the
properties known to hold in the $x$-independent case,
see~\cite{CP2020}.

\section{Peculiarities of the \texorpdfstring{$x$}{x}-Dependent Case}
\label{sec:pec}

The extension to the $x$-dependent case can not be merely reduced to
the rise of technical difficulties. Indeed, some properties are
irremediably lost and new phenomena arise, as shown below.

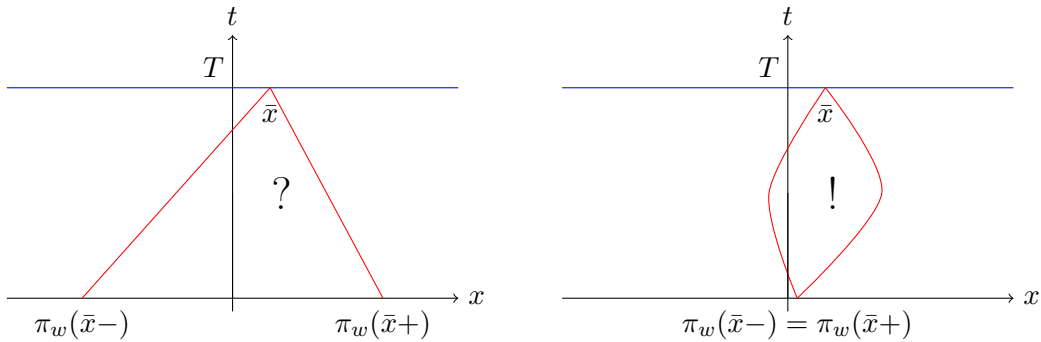
\begin{figure}[!ht]
  \begin{center}
    \begin{tikzpicture}[x=0.5cm,y=0.7cm]
      \draw[->] (-6, 0) -- (6, 0) node[right] {$x$};

      \draw[->] (0, -0.25) -- (0, 5) node[above] {$t$};

      \draw[-, blue] (-6, 4) -- (6, 4);

      \draw[-, red] (1, 4) -- (4,0);

      \draw[-, red] (1, 4) -- (-4,0);

      \node at (-4, -0.5) {$\pi_w (\bar x-)$};

      \node at (4, -0.5) {$\pi_w (\bar x+)$};

      \draw[-] (0,0) -- (0,2) node[right] {\quad\LARGE?};

      \node at (-0.5, 4.4) {$T$};

      \node at (1,3.5) {$\bar x$};

    \end{tikzpicture}
    \qquad
    \begin{tikzpicture}[x=0.5cm,y=0.7cm]

      \draw[->] (-6, 0) -- (6, 0) node[right] {$x$};

      \draw[->] (0, -0.25) -- (0, 5) node[above] {$t$};

      \draw[-, blue] (-6, 4) -- (6, 4);

      \draw[-, red] plot [smooth] coordinates {(1, 4) (-0.5, 2)
        (0.25,0)};

      \draw[-, red] plot [smooth] coordinates {(1, 4) (2.5, 2)
        (0.25,0)};

      \node at (0.25, -0.5) {$\pi_w (\bar x-)=\pi_w (\bar x+)$};

      \draw[-] (0,0) -- (0,2) node[right] {\quad\LARGE!};

      \node at (-0.5, 4.4) {$T$};

      \node at (1,3.5) {$\bar x$};

    \end{tikzpicture}
    \caption{Left, in the $x$-independent case, extremal
      characteristics are straight lines and those emanating from the
      point of jump $\bar x$ in $w$ at time $T$ select the segment
      $\mathopen]\pi_w (x-), \pi_w (x+)\mathclose[$ along the $x$ axis
      at time $0$ where the initial data has no effect on $w$. Right,
      in our $x$-dependent choice~\eqref{eq:11} of the flow,
      characteristics bend and uniquely determine the initial data
      evolving into $w$. Note that the solution in the region
      delimited by the characteristics is unique.}
    \label{fig:xtplane}
  \end{center}
\end{figure}
The most apparent difference between the two situations is described
in Figure~\ref{fig:xtplane}, with reference to extremal backward
generalized characteristics, whose behaviors in the two cases are
quite different. In the $x$-independent case, extremal backward
characteristics define a \emph{non uniqueness gap}, see
Figure~\ref{fig:xtplane}. On the contrary, in the $x$-dependent case,
extremal backward characteristics may well intersect at the initial
time, so that the non uniqueness gap disappears.

\begin{figure}[!h]
  \begin{center}
    \includegraphics[trim = {95 0 90 0}]{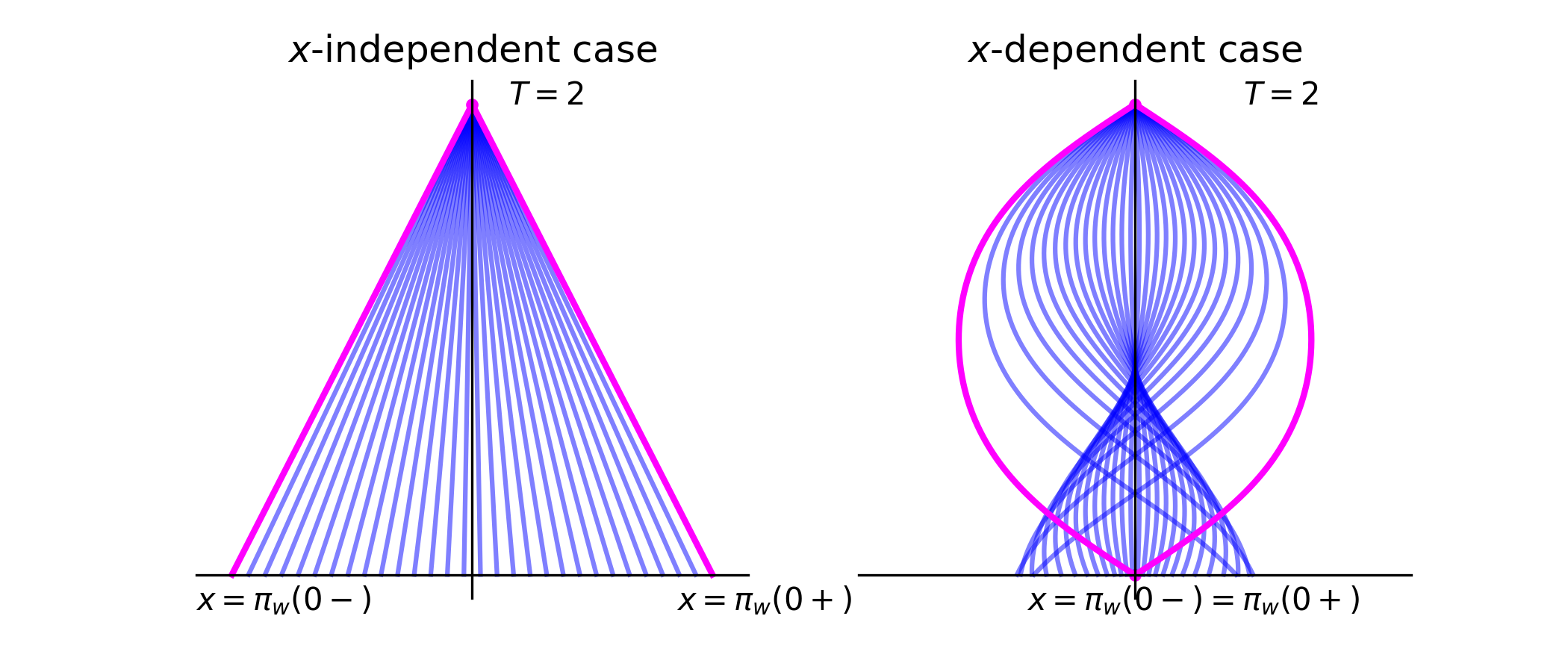}
  \end{center}
  \vspace*{-5mm}
  \caption{Left, in the $x$-independent case, the Hamiltonian rays
    fill the non uniqueness gap described in
    Figure~\ref{fig:xtplane}. Right, in the $x$-dependent case defined
    by the Hamiltonian~\eqref{eq:11}, extremal characteristics still
    do not intersect, but Hamiltonian rays do and may well exit the
    non uniqueness gap or also intersect.}
  \label{fig:shootingMethod}
\end{figure}

Furthermore, in the $x$-independent case, an isentropic solution,
see~\cite[Theorem~3.1]{CPS2022}, is constructed filling the non
uniqueness gap with Hamiltonian rays~\eqref{eq:hamiltonian_rays}
emanating from $q (T) = x$,
$p (T) = \theta \, w (x+) + (1-\theta) \, w (x-)$, for
$\theta \in [0,1]$. On the contrary, the same idea fails in the
$x$-dependent case. The numerical integrations in
Figure~\ref{fig:shootingMethod} referred to~\eqref{eq:ode_system} with
Hamiltonian~\eqref{eq:11}, show that extremal backward characteristics
still do not intersect in $\mathopen]0,T\mathclose[ \times \R$, but
the intermediate Hamiltonian rays may well cross each other and even
exit the region bounded by the extremal characteristics.

When $H$ does not depend on $x$, the $u_o^*$ defined in
Corollary~\ref{cor:geoCL} is characterized by~\cite[(G2) in
Proposition~5.2]{CP2020}. Then, \cite[(R1) in Lemma~7.2]{CP2020}
ensures not only that $u_o^*$ is one sided Lipschitz continuous, but
also that the solution $\tilde u$ to~\eqref{eq:cl} with datum $u_o^*$
evolving into $w$ is Lipschitz continuous on any compact subset of
$\mathopen]0, T \mathclose[ \times \R$. Thus, $\tilde u$ satisfies the
inequality in Definition~\ref{def:entropy_solution} with an equality,
i.e., it is an \emph{isentropic} and also reversible in time solution,
see related multi-dimensional results in~\cite{MR1722801}.

This actually characterizes the homogeneous case. Indeed, there exists
an $x$-dependent Hamiltonian $H$, a profile $w$ and a time $T > 0$
such that $I^{CL}_T(w) \neq \emptyset$ but in any solution evolving
from an initial datum in $I^{CL}_T(w)$ shocks arise at a time $t < T$,
so that no reversible solution is possible, see
Figure~\ref{fig:nice}. In other words, the profile $w$ can be reached
exclusively producing a sufficient amount of entropy and no isentropic
solution evolves into $w$. Each of these facts necessarily requires
$H$ to depend on $x$ and can not take place in an $x$-independent
setting, as shown in~\cite{CP2020}. A consequence is that no direct
definition of $u_o^*$ is available, as it was in the $x$-independent
case, and we have to resort to~\eqref{eq:hj} for its construction.

\begin{figure}[!h]
  \begin{center}
    \includegraphics[scale=0.45]{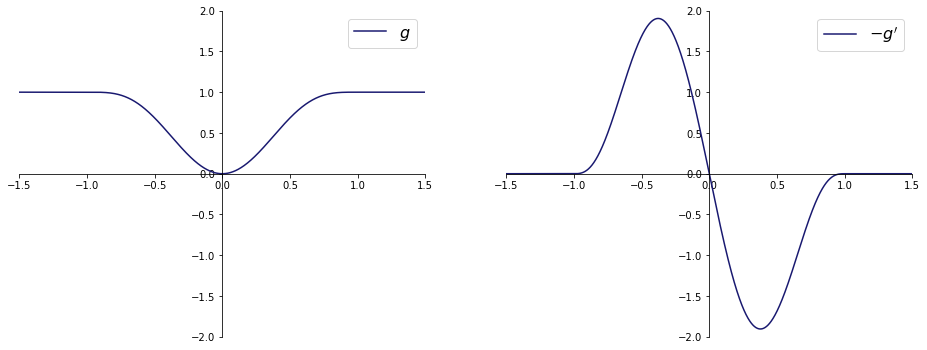}
  \end{center}
  \vspace*{-5mm}
  \caption{Left, graph of $g$ and, right, the graph of $-g'$,
    according to~\eqref{eq:11}. Clearly, $g$ is $\C3 (\R;[0,1])$,
    even, strictly increasing on $[0,1]$, $g'$ attains values in
    $[-2,2]$ and $H$ in~\eqref{eq:11} satisfies~\eqref{eq:space_dpcy}
    with $X=1$.}
  \label{fig:g}
\end{figure}

\begin{theorem}
  \label{th:counter_example}
  Define, see Figure~\ref{fig:g},
  \begin{equation}
    \label{eq:11}
    H(x, u) \coloneqq  \frac{u^2}{2} + g(x)
    \quad \mbox{ where } \quad
    g(x) \coloneqq
    \left\{
      \begin{array}{cl}
        1 - (1 - x^2)^4 & \text{if }\modulo{x} \leq 1 \,,
        \\
        1 & \text{otherwise.}
      \end{array}
    \right.
  \end{equation}
  Then, \eqref{eq:smoothness}, \eqref{eq:space_dpcy}
  and~\eqref{eq:convexity} hold. Moreover, for all
  $T > \pi / (2\sqrt2)$, there exists $w \in \L\infty (\R; \R)$ that
  contains a discontinuity and such that $I_T^{CL}(w)$ is a singleton.
\end{theorem}

\noindent The proof is deferred
to~\S~\ref{subsec:proof-theor-refth:c}.

Remark that if $H$ does not depend on $x$, as soon as $w$ has a jump,
then the contrary to the conclusion of
Theorem~\ref{th:counter_example} holds true,
see~\cite{CP2020}. Indeed, $I^{CL}_T (w)$ is either empty or infinite,
whenever $w$ has a discontinuity. In particular, \cite[(G1) in
Proposition~5.2]{CP2020} does not hold.

Recall that~\cite[Section~5]{MR1722801} presents, in the $n$
dimensional case, a backward procedure to construct what corresponds
here, in the $x$-independent case, to $U_o^*$
in~\eqref{eq:inverse_design}. Then, \cite[Example~6.3]{MR1722801}
proves that this procedure may well fail in the $x$-dependent case. In
Theorem~\ref{th:counter_example}, which is however restricted to the
$1$ dimensional case, the function $H$ also
satisfies~\eqref{eq:space_dpcy}, showing that the behavior for
$\modulo{x} \to +\infty$ is not relevant in this context. More
relevant, Theorem~\ref{th:counter_example} shows that there may well
be an \emph{intrinsic} minimal entropy production, independently of
any constructive procedure. As a matter of fact, the $U_o^*$
in~\eqref{eq:inverse_design} corresponds to the construction
in~\cite{MR1722801}, although it is built by means of optimal control
problems rather than by means of backward Hamilton--Jacobi
equations. However, we are here interested in the broader inverse
design characterization discussed in Section~\ref{sec:ext}, rather
than in time reversibility.

\begin{figure}[!ht]
  \begin{center}
    \includegraphics[scale=0.66]{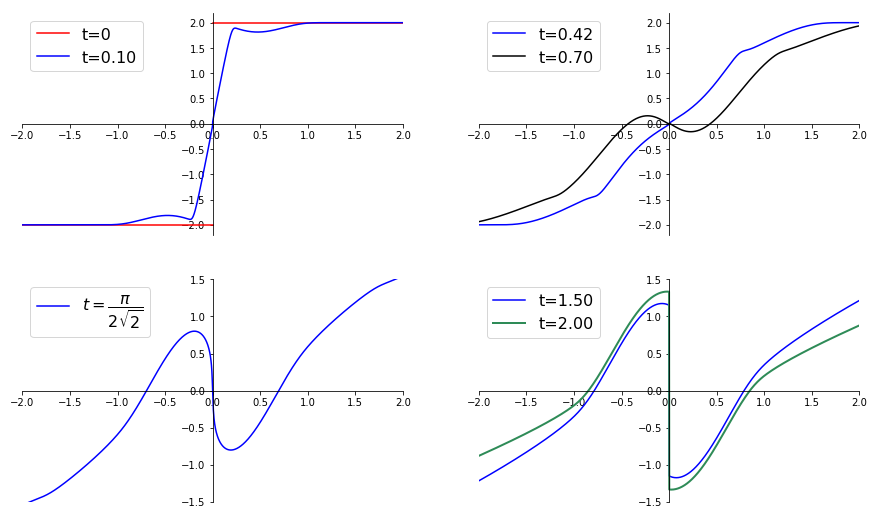}
  \end{center}
  \vspace*{-5mm}
  \caption{Evolution in time of (a numerical approximation of) the
    solution $u$ to~\eqref{eq:cl}--\eqref{eq:11}--\eqref{eq:12},
    constructed in Theorem~\ref{th:counter_example}, as a function of
    the space variable $x$, computed at different times, see also
    Figure~\ref{fig:nice}. Note the initial rarefaction profile
    turning into a shock at time
    $T = \left.\pi \middle/ (2\sqrt2)\right.$.}
  \label{fig:rarefaction_to_shock}
\end{figure}

The evolution of the numerical solution computed with a standard
finite volume scheme, is represented in
Figure~\ref{fig:rarefaction_to_shock}, see also
Figure~\ref{fig:nice}. Remark, and this is intrinsic to the
heterogeneous case, that the initial rarefaction profile evolves into
a shock wave.  The time asymptotic behavior shows further differences
with the $x$-independent case, see~\cite{procHYP2022} for more
details.

\section{Proofs}
\label{sec:inverse_design}

Several results of use below can be obtained through rather classical
techniques but can hardly be precisely localized in the literature. In
these cases, we refer to~\cite{SyllaPhD}, where all details are
provided.

\begin{proofof}{Lemma~\ref{lmm:regularity_profile_cl}}
  Let $u_o \in I^{CL} _T (w)$. Call $U_o$ a primitive of $u_o$, so
  that $\partial_x S^{HJ}_T U_o = w$
  by~\cite[Theorem~2.20]{CPS2022}. Then,
  by~\cite[Theorem~5.3.8]{CSBook}, $(t,x) \mapsto (S^{HJ}_t U_o) (x)$
  is locally semiconcave in the sense
  of~\cite[Definition~1.1.1]{CSBook}. Thus,
  $(t,x) \mapsto (S^{CL}_t u_o) (x)$ is locally one sided Lipschitz
  continuous in the space variable and hence in $\BVloc (\R;\R)$, for
  all $t>0$. As is well known, this ensures the existence of left and
  right traces at any point of the map $x \mapsto (S^{CL}_t u_o) (x)$,
  for all $t>0$.
\end{proofof}

\begin{proofof}{Proposition~\ref{prop:hypOk}}
  It is immediate to prove that~\eqref{eq:convexity}
  implies~\eqref{eq:WGNL}.  Thanks to~\eqref{eq:space_dpcy}
  and~\eqref{eq:convexity}, we can use~\cite[Lemma~8.1.3 and
  Corollary~8.1.4]{SyllaPhD} which ensure that there exists a function
  $\phi \in \C0(\R_+; \R)$ that verifies
  \begin{equation}
    \label{eq:10}
    \forall (x, p) \in \R^2, \quad H(x, p) \geq \phi\left(\modulo{p}\right)
    \quad \mbox{ with } \quad
    \frac{\phi(r)}{r} \limit{r}{+\infty} +\infty \,.
  \end{equation}
  Then, \eqref{eq:UC} readily follows.
\end{proofof}

\subsection{Proof of Theorem~\ref{th:two_min}}
\label{subsec:proof-theor-refth:tw}

\begin{lemma}
  \label{lmm:two_min}
  Let $H$ satisfy~\eqref{eq:smoothness}, \eqref{eq:space_dpcy}
  and~\eqref{eq:convexity}.  Fix $U_o \in \Lip(\R; \R)$ and let $U$
  solve~\eqref{eq:hj} in the sense of
  Definition~\ref{def:viscosity_solution}. Fix $T > 0$ and
  $\xi, \zeta \in \Lip([0, T]; \R)$, with $\xi \leq \zeta$. Then for
  all $s, \tau \in [0, T]$ with $s < \tau$,
  \begin{equation}
    \label{eq:two_min1}
    \begin{aligned}
      & \int_{\xi(\tau)}^{\zeta(\tau)} U(\tau, x) \d{x} -
      \int_{\xi(s)}^{\zeta(s)} U(s, x) \d{x} + \int_{s}^{\tau}
      \int_{\xi(t)}^{\zeta(t)} H \left(x, \partial_x U(t, x)\right)
      \d{x} \d{t}
      \\
      & = \int_{s}^{\tau} \left( \dot \zeta(t) \, U\left(t,
          \zeta(t)\right) - \dot \xi(t) \, U\left(t, \xi(t)\right)
      \right) \d{t}.
    \end{aligned}
  \end{equation}
\end{lemma}

\noindent This Lemma is analogous to~\cite[Lemma~3.2]{Dafermos1977},
see~\cite[Lemma~8.3.13]{SyllaPhD} for a detailed proof.

\begin{proofof}{Theorem~\ref{th:two_min}}
  We only prove the result for the maximal backward characteristic
  $\hat\zeta$, which we denote for simplicity $\zeta$. The details of
  the proof for the minimal characteristic $\check\zeta$ are similar.

  Fix $\epsilon > 0$. Apply Lemma~\ref{lmm:two_min} with $\zeta$ and
  $\xi = \zeta - \epsilon$ on $[0, t]$ and $s=0$, $\tau = t$.  After
  dividing by $\epsilon$, we obtain:
  \begin{equation}
    \label{eq:two_min2}
    \begin{aligned}
      & \frac{1}{\epsilon} \int_{\zeta(t) - \epsilon}^{\zeta(t)} U(t,
      y) \d{y} - \frac{1}{\epsilon} \int_{\zeta(0) -
        \epsilon}^{\zeta(0)} U_o(y) \d{y} + \frac{1}{\epsilon}
      \int_0^{t} \int_{\zeta(s) - \epsilon}^{\zeta(s)} H\left(y,
        \partial_x U(s, y)\right) \d{y} \d{s}
      \\
      & = \frac{1}{\epsilon} \int_0^{t} \dot \zeta(s) \left( U\left(s,
          \zeta(s)\right) - U\left(s, \zeta(s) - \epsilon\right)
      \right) \d{s}.
    \end{aligned}
  \end{equation}
  We want to pass to the limit $\epsilon \to 0$
  in~\eqref{eq:two_min2}. To this aim, recall that $U (t,\cdot)$ and
  $U_o$ are continuous in $x$ by
  Definition~\ref{def:viscosity_solution}. Moreover, $\partial_x U$
  solves~\eqref{eq:cl} in the sense of
  Definition~\ref{def:entropy_solution} with initial data $U_o'$,
  see~\cite[Theorem~2.20]{CPS2022}. For a.e.~$s \in [0, t]$,
  $\partial_x U(s,\cdot)$ has left and right limits at $x = \zeta(s)$
  that exist by Lemma~\ref{lmm:regularity_profile_cl} and coincide,
  since $\zeta$ is genuine~\cite[Definition~3.2 and
  Theorem~3.2]{Dafermos1977}. The map $U$ is Lipschitz continuous,
  hence
  $U (s,x_2) - U (s,x_1) = \int_{x_1}^{x_2} \partial_x U (s,y) \,
  \d{y}$, so that for a.e.~$s \in [0, t]$, $U(s, \cdot)$ is
  differentiable at $x = \zeta(s)$ and
  $\partial_x U\left(s, \zeta(s)\right) = \lim_{h \to 0} \partial_x
  U(s, \zeta(s) + h)$.  We thus obtain:
  \begin{equation}
    \label{eq:two_min3}
    U\left(t, \zeta(t)\right) - U_o\left(\zeta(0)\right)
    = \int_0^t \left(
      \dot \zeta(s) \, \partial_x U\left(s, \zeta(s)\right)
      -
      H\left(
        \zeta(s), \partial_x U\left(s, \zeta(s)\right)
      \right) \right) \d{s}.
  \end{equation}
  Since $\zeta$ is genuine, there exists a function
  $\omega \in \C{1}(\mathopen]0, t\mathclose[; \R)$ such that
  $(\zeta, \omega)$ is a solution to system~\eqref{eq:ode_system} with
  final conditions $\zeta(t) = x$ and
  $\omega(t) = \partial_x U(t, x+)$, since $\zeta$ is a maximal
  characteristics, see~\cite[Theorem~3.3]{Dafermos1977}. Moreover, for
  a.e.~$s \in [0, t]$,
  $\omega(s) = \partial_x U\left(s, \zeta(s)\right)$. Combining these
  details with~\eqref{eq:two_min3} and~\eqref{eq:convexity}, classical
  computations lead to:
  \begin{displaymath}
    \begin{aligned}
      U(t, x) - U_o\left(\zeta(0)\right) & = \int_0^t \left( \dot
        \zeta(s) \; \omega(s) - H\left(\zeta(s), \omega(s) \right)
      \right) \d{s}
      \\
      & = \int_0^t \left( \partial_p H \left(\xi(s), \omega(s)\right)
        \, \omega(s) - H\left(\zeta(s), \omega(s) \right) \right)
      \d{s}
      \\
      & = \int_0^t L\left(\zeta(s), \partial_p H \left(\xi(s),
          \omega(s)\right) \right) \d{s}
      \\
      & = \int_0^t L\left(\zeta(s), \dot \zeta(s) \right) \d{s} \,,
    \end{aligned}
  \end{displaymath}
  concluding the proof.
\end{proofof}

\noindent Further information about the regularity of $U$ along
characteristics can be found in~\cite[\S~5.5]{CSBook}.

\subsection{Proof of Theorem~\ref{th:nonemptiness}}
\label{subsec:proofs-relat-sect}

In the light of the regularity proved in
Lemma~\ref{lmm:regularity_profile_cl}, whenever
$w \in \L{\infty}(\R; \R)$ is such that $I_T^{CL}(w) \neq \emptyset$,
then by $w(x)$, we mean the left trace $w (x-)$ of $w$ at $x$, for all
$x \in \R$.

\begin{lemma}
  \label{lem:HS}
  Let $H$ satisfy~\eqref{eq:smoothness}, \eqref{eq:space_dpcy}
  and~\eqref{eq:convexity}. Then, for all $(q_o,p_o) \in \R^2$, the
  Cauchy problem~\eqref{eq:ode_system} with initial datum $(q_o,p_o)$
  at time $0$ admits a unique maximal solution $(q,p)$ defined on all
  $\R$ and satisfying, with the notation~\eqref{eq:legendre},
  \begin{equation}
    \label{eq:3}
    \modulo{p (t)}
    \leq
    \sup_{\substack{\modulo{u} \leq \modulo{p_o} \\ x \in \R}} \modulo{H(x, u)}
    + \sup_{\substack{\modulo{v} \leq 1 \\ x \in \R}} L(x, v) \,.
  \end{equation}
  Moreover, calling $(q,p)$ the solution to~\eqref{eq:ode_system} with
  datum $(q_o,p_o)$ at time $0$, the maps
  \begin{equation}
    \label{eq:flow}
    \!\!\!\!\!\!\!
    \fonction{\flow}{\R^3}{\R^2}{(t, q_o, p_o)}{\left(q(t), p(t)\right)}\;
    \fonction{\flow_q}{\R^3}{\R}{(t, q_o, p_o)}{q(t)}\;
    \fonction{\flow_p}{\R^3}{\R}{(t, q_o, p_o)}{p(t)}
    \!\!
  \end{equation}
  are of class $\C2$.
\end{lemma}

\begin{proofof}{Lemma~\ref{lem:HS}}
  By~\eqref{eq:smoothness}, the standard Cauchy Lipschitz Theorem
  ensures local existence and uniqueness of a solution $(q,p)$ to the
  Cauchy problem for~\eqref{eq:ode_system} with datum
  $(q_o,p_o)$. Moreover, since $H$ is conserved along solutions
  to~\eqref{eq:ode_system}, for all $t$ where $(q,p)$ is defined,
  \begin{displaymath}
    \sup_{\substack{\modulo{u} \leq \modulo{p_o} \\ x \in \R}} \modulo{H(x, u)}
    \geq H(q_o, p_o)
    = H\left(q(t), p(t)\right)
    \geq \modulo{p(t)}
    - \sup_{\substack{\modulo{v} \leq 1 \\ x \in \R}} L(x, v) \,,
  \end{displaymath} where we used~\eqref{eq:legendre}, see also~\cite[Formula~(8.1.5)]{SyllaPhD}
  with $\lambda=1$, proving~\eqref{eq:3}. By~\eqref{eq:ode_system},
  \eqref{eq:smoothness} and~\eqref{eq:3}, we also have that the solution
  $(q,p)$ is bounded and uniformly continuous on bounded
  intervals. Hence, it is globally defined.

  Standard results on ordinary differential equations, see
  e.g.~\cite[Theorem~3.9, Theorem~3.10]{MR2347697}, ensure that the
  flow $\flow$ is as regular as $\partial_q H$, $\partial_x H$ and,
  by~\eqref{eq:smoothness}, the proof is completed.
\end{proofof}

The next three lemmas state in full rigor simple geometric properties
that are consequences of~\eqref{eq:convexity}
and~\eqref{eq:space_dpcy} on the graph of $H$ (essentially, a
\emph{canyon} along the $x$ direction).

\begin{lemma}
  \label{lem:u}
  Let $H$ satisfy~\eqref{eq:smoothness}, \eqref{eq:space_dpcy}
  and~\eqref{eq:convexity}.  Then, there exists a unique function
  \begin{equation}
    \label{eq:sm1}
    \fonction{\u}{\R}{\R}{x}{\u (x)}
    \quad \mbox{ such that } \quad
    \forall\, x \in \R \quad
    \partial_p H\left(x, \u(x)\right) = 0 \,.
  \end{equation}
  Moreover, $\u \in \C2 (\R; \R)$, if $\modulo{x} \geq X$ then
  $\u' (x) = 0$ and the following quantities are well defined
  \begin{equation}
    \label{eq:sm2}
    \underline{\u} \coloneqq  \min_{x \in \R} \u(x)
    \,, \qquad
    \overline{\u} \coloneqq  \max_{x \in \R} \u(x)
    \,, \qquad
    \bK \coloneqq  \max_{x \in \R} H\left(x, \u(x)\right) \,.
  \end{equation}
\end{lemma}

\begin{proofof}{Lemma~\ref{lem:u}}
  Existence and uniqueness of $\u$ follow
  from~\eqref{eq:convexity}. Together, \eqref{eq:smoothness}
  and~\eqref{eq:convexity} allow to apply the Implicit Function
  Theorem, proving both the $\C2$ regularity of $\u$ and,
  by~\eqref{eq:space_dpcy}, that $\u' (x) = 0$ whenever
  $\modulo{x} \geq X$. The completion of the proof is now immediate.
\end{proofof}

\begin{lemma}
  \label{pp:sm1}
  Let $H$ satisfy~\eqref{eq:smoothness}, \eqref{eq:space_dpcy}
  and~\eqref{eq:convexity}.  Referring to the function $\u$ and to the
  constant $\bK$ defined in Lemma~\ref{lem:u},
  there exist functions
  \begin{equation}
    \label{eq:6}
    \fonction{m}{\R \times \mathopen]\bK, +\infty\mathclose[}{\R}{(x,c)}{m (x,c)}
    \quad \mbox{ and } \quad
    \fonction{M}{\R \times \mathopen]\bK, +\infty\mathclose[}{\R}{(x,c)}{M (x,c)}
  \end{equation}
  uniquely characterized, for $c > \bK$ and $x \in \R$, by
  \begin{equation}
    \begin{array}{c}
      H\left(x, m (x,c)\right) = c \quad \mbox{ and } \quad m (x,c) < z (x)
      \\
      H\left(x, M (x,c)\right) = c \quad \mbox{ and } \quad M (x,c) > z (x)
    \end{array}
    \label{eq:7}
  \end{equation}
  Moreover,
  \begin{enumerate}[\bf(i)]
  \item\label{item:4}
    $m, M \in \C{1}(\R \times \mathopen]\bK, +\infty\mathclose[;\R)$.

  \item \label{item:5} $m$ and $M$ have a compact space dependency:
    \begin{equation}
      \label{eq:sm11}
      \modulo{x} > X \mbox{ and } c > \bK
      \implies
      \partial_x m (x, c) = 0
      \text{ and }
      \partial_x M (x, c) = 0 \,.
    \end{equation}

  \item\label{item:6} For all $x \in \R$, $m(x, \cdot)$ is decreasing
    while $M(x, \cdot)$ is increasing.

  \item\label{item:7} For all $x \in \R$,
    $\lim_{c\to+\infty} m (x,c) = -\infty$ and
    $\lim_{c\to+\infty} M (x,c) = +\infty$.
  \end{enumerate}
\end{lemma}

\begin{proofof}{Lemma~\ref{pp:sm1}}
  We only prove the results for $M$, the details for $m$ are entirely
  similar.

  Assumption~\eqref{eq:convexity} ensures that condition~\eqref{eq:7}
  uniquely defines the map $M$ in~\eqref{eq:6}. An application of the
  Implicit Function Theorem shows the regularity,
  by~\eqref{eq:smoothness}, proving~\eqref{item:4}.

  Again, by the Implicit Function Theorem and the chain rule, we have
  \begin{displaymath}
    \partial_x H \left(x, M(x, c)\right)
    +
    \partial_x M(x, c) \; \partial_p H\left(x, M(x, c)\right)
    = 0 \,,
  \end{displaymath}
  which implies~\eqref{eq:sm11} by~\eqref{eq:space_dpcy} and
  by~\eqref{eq:6} and~\eqref{eq:7}, proving~\eqref{item:5}.

  By the definitions~\eqref{eq:sm1} of $\u$
  and~\eqref{eq:6}--\eqref{eq:7} of $M$, we have that for all
  $(x, c) \in \R \times \mathopen]\bK, +\infty\mathclose[$,
  $M(x, c) > \u(x)$ and hence
  $\partial_p H\left(x, M(x, c)\right) > 0$. Again using the Implicit
  Function Theorem,
  $\partial_c M(x, c) \; \partial_p H\left(x, M(x, c)\right) = 1$
  proving that $\partial_c M(x, c)>0$, proving~\eqref{item:6}.

  Since $M$ is increasing,
  $\lim_{c\to+\infty} M (x,c) = \sup_{c>\bK} M (x,c)$. The boundedness
  of $c \mapsto M (x,c)$ for some $x$ then contradicts the equality
  $H\left(x, M(x, c)\right) = c$, proving~\eqref{item:7}.
\end{proofof}

\begin{lemma}
  \label{pp:sm2}
  Let $H$ satisfy~\eqref{eq:smoothness}, \eqref{eq:space_dpcy}
  and~\eqref{eq:convexity}.  Referring to the constant $\bK$ in
  Lemma~\ref{lem:u} and to the functions $m,M$ in Lemma~\ref{pp:sm1},
  define the functions:
  \begin{displaymath}
    \fonction{v}{\mathopen]\bK, +\infty\mathclose[}{\R}{c}{\sup_{x \in \R} \partial_p H\left(x, m(x, c)\right)}
    \; \mbox{ and }\;
    \fonction{V}{\mathopen]\bK, +\infty\mathclose[}{\R}{c}{\inf_{x \in \R} \partial_p H\left(x, M(x, c)\right)}
  \end{displaymath}
  Then:
  \begin{enumerate}[\bf(i)]
  \item\label{item:8} $v$ is nonincreasing and $V$ is nondecreasing;

  \item\label{item:9} $\lim_{c\to+\infty} v (c) = -\infty$ and
    $\lim_{c\to+\infty} V (c) = +\infty$.
  \end{enumerate}
\end{lemma}

\begin{proofof}{Lemma~\ref{pp:sm2}}
  By~\eqref{eq:space_dpcy} and~\eqref{item:5} in
  Lemma~\ref{pp:sm1}, 
  $v$ and $V$ are well-defined.  We now prove the
  statements~\eqref{item:8} and~\eqref{item:9} for $V$, the case of
  $v$ being entirely analogous.

  From the monotonicity of $M$ and $\partial_p H$ with respect to
  their second argument:
  \begin{displaymath}
    \forall \, x \in \R, \; \forall \, c_1, c_2 > \bK ,\;
    \mbox{ if } c_1 < c_2 \mbox{ then }
    \partial_p H\left(x, M(x, c_1)\right)
    \leq
    \partial_p H\left(x, M(x, c_2)\right).
  \end{displaymath}
  Taking the infimum over $x \in \R$ we prove~\eqref{item:8}.

  By~\eqref{item:8}, $\lim_{c\to+\infty} V (c) = \sup_{c > \bK}
  V(c)$. By contradiction, assume that
  $\overline{V}\coloneqq \sup_{c > \bK} V(c)$ is finite. By the
  definition of $V$ and~\eqref{eq:space_dpcy},
  \begin{equation}
    \label{eq:sm4}
    \forall \, n \in \N  \cap \mathopen]\bK, +\infty\mathclose[
    \quad
    \exists \,x_n \in [-X, X],
    \quad \mbox{ such that } \quad
    \partial_p H\left(x_n, M(x_n, n)\right) \leq \overline{V} \,.
  \end{equation}
  Up to a subsequence, we can assume that $(x_n)_n$ converges to some
  $\overline{x} \in [-X, X]$. Item~\eqref{item:7} in
  Lemma~\ref{pp:sm1} and~\eqref{eq:convexity} imply that
  $\lim_{n \to +\infty} \partial_p H\left(\overline{x},
    M(\overline{x}, n)\right) = +\infty$. Therefore,
  \begin{equation}
    \label{eq:sm5}
    \exists \, N \in \N \cap \mathopen]\bK, +\infty\mathclose[
    \quad \mbox{ such that } \quad
    \forall \, n \geq N, \quad
    \partial_p H\left(\overline{x}, M(\overline{x}, n)\right)
    >
    \overline{V} \,.
  \end{equation}
  The monotonicity of $M$ and~\eqref{eq:convexity}, combined
  with~\eqref{eq:sm4}, result in:
  \begin{displaymath}
    \forall n \geq N, \quad \partial_p H\left(x_n, M(x_n, N)\right)
    \leq \overline{V}
  \end{displaymath}
  which contradicts~\eqref{eq:sm5}, proving~\eqref{item:9}.
\end{proofof}

\begin{lemma}
  \label{lmm:surjectivity_sf}
  Let $H$ satisfy~\eqref{eq:smoothness}, \eqref{eq:space_dpcy}
  and~\eqref{eq:convexity}.  Fix $q_o \in \R$ and $T \in \R$. Then,
  the map $p \mapsto \flow_q (T,q_o,p)$ defined in~\eqref{eq:flow} is
  surjective, in the sense that
  \begin{equation}
    \label{eq:surjectivity_sf}
    \forall \, q_T \in \R, \quad
    \exists \, p_* \in \R \quad
    \mbox{ such that } \quad
    \flow_q(T, q_o, p_*) = q_T \,,
  \end{equation}
  or, with the notation~\eqref{eq:hamiltonian_rays},
  \begin{equation}
    \label{eq:9}
    \forall \, (x_o, x_T) \in \R^2, \quad
    \exists \, q \in \mathcal{R}_T \quad
    \mbox{ such that } \quad
    q (0) = x_o \mbox{ and } q (T) = x_T \,.
  \end{equation}
\end{lemma}

\begin{proofof}{Lemma~\ref{lmm:surjectivity_sf}}
  Recall the map $\u$ and the scalar $\bK$ defined in
  Lemma~\ref{lem:u}. Fix $c_o> \bK$ and let $p_o = M (q_o,c_o)$. By
  the conservation of $H$, for all $t \in \R$ we have
  $H\left(\flow_q(t, q_o, p_o), \flow_p(t, q_o, p_o)\right) = c_o$.
  Hence, by the definition~\eqref{eq:sm2} of $\bK$, for all
  $t \in \R$, it follows that
  $\flow_p(t, q_o, p_o) \neq \u\left(\flow_q(t, q_o, p_o)\right)$.

  Using the continuity of $\flow$, proved in Lemma~\ref{lem:HS}, as
  well as the fact that $p_o > \u(q_o)$, we deduce that
  \begin{displaymath}
    \forall t \in \R, \quad
    \flow_p(t, q_o, p_o) > \u\left(\flow_q(t, q_o, p_o)\right) \,.
  \end{displaymath}
  Therefore, for all $t \in \R$,
  $\flow_p(t, q_o, p_o) = M\left(\flow_q(t, q_o, p_o),c_o\right)$, as
  defined in~\eqref{eq:6}--\eqref{eq:7}. Thus,
  by~\eqref{eq:ode_system} and the definition of $V$ in
  Lemma~\ref{pp:sm2}, we have:
  \begin{displaymath}
    \flow_q(T, q_o, p_o)
    =
    q_o + \int_0^T \dot q(s) \d{s}
    \geq
    q_o +  V(c_o) \, T
    \limit{c_o}{+\infty}
    +\infty \,.
  \end{displaymath}
  A similar argument, with $p_o = m (q_o,c_o)$, yields
  $\flow_q(T, q_o, p_o) \leq q_o + v(c_o) \, T \limit{c_o}{+\infty}
  -\infty$. The continuity of $\flow_q$ coupled with the Intermediate
  Value Theorem concludes the proof of
  Lemma~\ref{lmm:surjectivity_sf}.
\end{proofof}

\begin{lemma}
  \label{lmm:candidat_inverse_design2}
  Let $H$ satisfy~\eqref{eq:smoothness}, \eqref{eq:space_dpcy}
  and~\eqref{eq:convexity}. Fix $T>0$ and $W \in \Lip(\R; \R)$ such
  that $I_T^{HJ}(W) \neq \emptyset$. Then, for all
  $U_o \in I_T^{HJ}(W)$, with the notation~\eqref{eq:inverse_design},
  \begin{equation}
    \label{eq:candidat_inverse_design2}
    \forall x \in \R, \quad U_o(x) \geq U_o^{*}(x).
  \end{equation}
\end{lemma}

\begin{proofof}{Lemma~\ref{lmm:candidat_inverse_design2}}
  Fix $x \in \R$ and $y \in \cR_T$ so that $y(0) = x$. Since
  $U_o \in I_T^{HJ}(W)$, by~\eqref{eq:24} we have:
  \begin{eqnarray*}
    &
    & W\left(y(T)\right) - \int_0^T L\left(y(s), \dot y(s)\right)
      \d{s}
    \\
    & =
    & \inf_{\substack{\gamma(T) = y(T) \\ \gamma \in \cR_T}}
    \left( \int_0^T L\left(\gamma(s), \dot \gamma(s)\right) \d{s} +
    U_o\left(\gamma(0)\right) \right) - \int_0^T L\left(y(s), \dot
    y(s)\right) \d{s}
    \\
    & \leq
    & \int_0^T L\left(y(s), \dot y(s)\right) \d{s} +
      U_o\left(y(0)\right) - \int_0^T L\left(y(s), \dot y(s)\right)
      \d{s}
    \\
    & \leq
    & U_o\left(y(0)\right)
    \\
    & =
    & U_o(x) \,.
  \end{eqnarray*}
  By taking the supremum over $y \in \cR_T$,
  by~\eqref{eq:inverse_design} we complete the proof.
\end{proofof}

By the second condition in~\eqref{eq:10}, for any
$W \in \Lip(\R; \R)$, there exists $\C{}_{H, W} > 0$ such that
\begin{equation}
  \label{eq:the_constant}
  \forall r \in \R_+, \; r \geq \C{}_{H, W} \implies \frac{\phi(r)}{1 + r} >
  \biggl( \sup_{\substack{\modulo{p} \leq \norma{W'}_{\L{\infty}(\R; \R)} \\ q \in \R}}
  \modulo{H(q, p)}
  + \sup_{\substack{\modulo{v} \leq 1 \\ q \in \R}} \modulo{L(q, v)} \biggr),
\end{equation}

\begin{lemma}
  \label{lmm:cs2}
  Let $H$ satisfy~\eqref{eq:smoothness}, \eqref{eq:space_dpcy}
  and~\eqref{eq:convexity}. Fix $T>0$ and $W \in \Lip(\R; \R)$. Then,
  $U_o^{*}$ defined by~\eqref{eq:inverse_design} is Lipschitz
  continuous and
  \begin{equation}
    \label{eq:candidat_inverse_design1}
    \norma{(U_o^{*})'}_{\L{\infty}(\R; \R)}
    \leq T \biggl( \sup_{\substack{\modulo{v} \leq \C{}_{H, W} \\ x \in \R}}
    \modulo{\partial_x L (x, v)} + \norma{W'}_{\L{\infty}(\R; \R)} \biggr) \,.
  \end{equation}
  Moreover, in its definition~\eqref{eq:inverse_design}, the $\sup$ is
  attained and for any Hamiltonian ray $q$ realizing the maximum
  in~\eqref{eq:inverse_design},
  $\norma{\dot q}_{\L{\infty}([0,T];\R)} \leq \C{}_{H, W}$.
\end{lemma}

\begin{proofof}{Lemma~\ref{lmm:cs2}}
  First, thanks to~\eqref{eq:inverse_design}, remark that
  \begin{displaymath}
    \forall x \in \R, \quad
    - U_o^{*}(x) =
    \inf_{\substack{q(0) = x \\ q \in \cR_T}}
    \left(
      \int_0^T L\left(q(s), \dot q(s)\right) \d{s} - W\left(q(T)\right)
    \right).
  \end{displaymath}
  Now reverse time applying the change of variable
  $\tau \coloneqq T - s$ and introducing
  $q^r (\tau) \coloneqq q (T-\tau)$,
  $L^r (x, v) \coloneqq L (x, - v)$,
  $H^r (x,p) \coloneqq \sup_{v \in \R} \left(p\, v\ - L^r (x,v)
  \right)$ so that $H^r (x,p) = H (x, -p)$. Moreover, using
  $p^r (\tau) = -p (T-\tau)$, $q \in \mathcal{R}_T$ if and only if
  $q^r \in \mathcal{R}^r_T$, where $\mathcal{R}^r_T$ is the set of
  Hamiltonian rays~\eqref{eq:hamiltonian_rays} defined by $H^r$, with
  reversed time.
  \begin{displaymath}
    \forall x \in \R, \quad
    - U_o^{*}(x) = \inf_{\substack{q^r(T) = x \\ q^r \in \cR^r_T}}
    \left( \int_0^T L^r\left(q^r(\tau), \dot q^r(\tau)\right) \d{\tau} -
      W\left(q^r(0)\right) \right) \,.
  \end{displaymath}
  Then, in view of~\cite[Corollary~8.3.15]{SyllaPhD},
  \begin{equation}
    \label{eq:25}
    \forall x \in \R, \quad
    - U_o^{*}(x) = \inf_{\substack{q^r(T) = x \\ q^r \in \Lip ([0,T];\R)}}
    \left( \int_0^T L^r\left(q^r(\tau), \dot q^r(\tau)\right) \d{\tau} -
      W\left(q^r(0)\right) \right) \,.
  \end{equation}
  So that, by~\cite[Theorem~8.3.12]{SyllaPhD},
  $-U_o^* (x) = U^r (T,x)$, with $U^r$ being the viscosity solution to
  the Hamilton-Jacobi equation
  \begin{displaymath}
    \left\{
      \begin{array}{l}
        \partial_t U^r + H^r (x, \partial_x U^r) = 0
        \\
        \tilde U (0,x) = -W (x) \,,
      \end{array}
    \right.
  \end{displaymath}
  proving~\eqref{eq:candidat_inverse_design1}.

  The result in~\cite[Corollary~8.3.15]{SyllaPhD} ensures that the
  supremum in Definition~\eqref{eq:inverse_design} is attained as a
  maximum. We can now combine~\cite[Theorem~8.3.9]{SyllaPhD}
  and~\cite[Corollary~8.3.15]{SyllaPhD} to complete the proof.
\end{proofof}

Remarkably, the next Lemma does not require $W$ to be reachable.

\begin{lemma}
  \label{lmm:cs3}
  Let $H$ satisfy~\eqref{eq:smoothness}, \eqref{eq:space_dpcy}
  and~\eqref{eq:convexity}. Fix $T>0$ and $W \in \Lip(\R; \R)$. Then
  $\cG$, as defined in~\eqref{eq:graph} has the following properties.
  \begin{enumerate}[(i)]
  \item\label{item:10} $\cG$ is surjective in the following sense:
    \begin{equation}
      \label{eq:graph1}
      \forall x_o \in \R, \; \exists x_T \in \R, \quad (x_o, x_T) \in \cG.
    \end{equation}
  \item\label{item:11} $\cG$ is a closed subset of $\R^2$.
  \item\label{item:12} For all $(x_o, x_T) \in \cG$, we have
    \begin{equation}
      \label{eq:graph2}
      \modulo{x_o - x_T} \leq T \, \C{}_{H, W} \,.
    \end{equation}
  \item\label{item:13} $\cG$ is monotone in the following sense:
    \begin{equation}
      \label{eq:graph_monotone}
      \forall \, (x_o, x_T), (y_o, y_T) \in \cG, \quad
      \begin{array}{r@{\,}c@{\,}l@{\implies}r@{\,}c@{\,}l}
        x_o
        & <
        & y_o
        & x_T
        & \leq
        & y_T\,;
        \\
        x_T
        & <
        & y_T
        & x_o
        & \leq
        & y_o \,.
      \end{array}
    \end{equation}
  \end{enumerate}
\end{lemma}

\noindent It is worth noting the connection between~\eqref{item:13}
and~\cite[Lemma~8.1]{CP2020}, inspired by~\cite[Lemma~2.2]{Lax}, see
also~\cite[Sections~5 and~6]{MR47234}. However, the $x$ dependence
makes the present proof significantly different and more intricate.

\begin{proofof}{Lemma~\ref{lmm:cs3}}
  Consider the different items separately.

  Property~\eqref{eq:graph1} comes from the definition of $U_o^{*}$,
  since the $\sup$ is actually a $\max$ by Lemma~\ref{lmm:cs2}.

  \paragraph{Proof of~\eqref{item:11}:} Let $(x_o^{n}, x_T^{n})_n$ be
  a sequence taking values in $\cG$ which converges to some
  $(x_o, x_T) \in \R^2$. By definition, for all $n \in \N$, there
  exists $q^{n} \in \cR_T$ such that
  \begin{equation}
    \label{eq1:inverse_design3}
    x_o^{n} = q^{n}(0), \quad x_T^{n} = q^{n}(T), \quad
    U_o^{*}(x_o^{n}) = W(x_T^{n})
    - \int_0^T L\left(q^{n}(s), \dot q^{n}(s)\right) \d{s} \,.
  \end{equation}
  For all $n \in \N$, let us denote by $p^{n} \in \C{1}([0, T];\R)$ a
  curve associated with $q^{n}$, given by
  $(x_o^{n}, x_T^{n}) \in \cG$.  Lemma~\ref{lmm:cs2} ensures that for
  all $n \in \N$,
  $\norma{\dot q^{n}}_{\L{\infty}(]0,T[;\R)} \leq \C{}_{H, W}$.  Note
  that for all $n \in \N$, $q^n \in \C1 ([0,T];\R)$
  by~\eqref{eq:hamiltonian_rays}. Thanks to~\eqref{eq:convexity}
  and~\eqref{eq:legendre},
  \begin{displaymath}
    \dot q^n (t) = \partial_p H\left(q^{n}(t), p^{n}(t)\right)
    \iff
    p^n (t) = \partial_v L\left(q^n (t),\dot q^n (t)\right) \,.
  \end{displaymath}
  This proves the boundedness of $\left(p^n (0)\right)$ and, up to a
  subsequence, we can assume that $\left(p^{n}(0), q^{n}(0)\right)_n$
  converges to $(p_o, x_o)$ with $p_o \in \R$. By Lemma~\ref{lem:HS},
  the flow of the Hamiltonian system is continuous and we establish
  the existence of $(q, p) \in \C0([0, T];\R^2)$ such that $(q^{n})_n$
  and $(p^{n})_n$ converge uniformly on $[0, T]$ to $q$ and $p$,
  respectively. Using the integral form of~\eqref{eq:ode_system}, we
  deduce that $(q,p)$ solves~\eqref{eq:ode_system}. Hence,
  $q \in \cR_T$ and $(x_o, x_T) \in \cG$.

  \paragraph{Proof of~\eqref{item:12}:} Let $(x_o, x_T) \in \cG$ and
  let $q \in \cR_T$ with $q (0) = x_o$ and $q (T) = x_T$.  Then, in
  view of Lemma~\ref{lmm:cs2} (latter part), we have
  \begin{displaymath}
    \modulo{x_o - x_T}
    = \modulo{q(0) - q(T)}
    \leq T \, \norma{\dot q}_{\L{\infty}(]0,T[;\R)}
    \leq T \, \C{}_{H, W} \,.
  \end{displaymath}

  \paragraph{Proof of~\eqref{item:13}:} We only prove the first
  implication in~\eqref{eq:graph_monotone}, the details of the proof
  for the second one are similar so we omit them.  Let
  $x, y \in \cR_T$ be two maximizers for $U_o^{*}(x_o)$ and
  $U_o^{*}(y_o)$, respectively. By assumption, we have $x(0) < y(0)$
  so that we can define
  \begin{equation}
    \label{eq:monotone1}
    \tau
    =
    \sup \left\{
      t \in [0,T] \colon x (s) < y (s) \mbox{ for all }s \in [0,t]
    \right\}
  \end{equation}
  and assume, by contradiction, that $\tau < T$, so that
  $x (\tau) = y (\tau)$.  Define the concatenation
  \begin{displaymath}
    \forall t \in [0, T], \quad
    \xi(t) =
    \left\{
      \begin{array}{lcl}
        y(t) & \text{if} & 0 \leq t \leq \tau\,,
        \\
        x(t) & \text{if} & \tau < t \leq T \,.
      \end{array}
    \right.
  \end{displaymath}
  Clearly, $\xi \in \Lip([0,T];\R)$ and $\xi(0) = y_o$.

  We now prove that $\dot x(\tau) \neq \dot y(\tau)$. Denote
  $p_x, p_y \in \C{1}([0, T]; \R)$ the curves associated with $x$ and
  $y$, respectively, given by $(x_o, x_T), (y_o, y_T) \in \cG$. Then,
  \begin{displaymath}
    \begin{aligned}
      \dot x(\tau) = \dot y(\tau) & \iff \partial_p H\left(x(\tau),
        p_x(\tau)\right) = \partial_p H\left(y(\tau), p_y(\tau)\right)
      \\
      & \iff \partial_p H\left(y(\tau), p_x(\tau)\right) = \partial_p
      H\left(y(\tau), p_y(\tau)\right) \iff p_x(\tau) = p_y(\tau) \,,
    \end{aligned}
  \end{displaymath}
  since $p \mapsto \partial_p H\left(y(\tau), p\right)$ is a bijection
  by~\eqref{eq:convexity}. However this contradicts the uniqueness of
  solutions to~\eqref{eq:ode_system}, see Lemma~\ref{lem:HS}.

  Hence, $\xi$ is not differentiable at point $\tau$. Moreover, since
  $x$ and $y$ are maximizers, we have, in light of the Dynamic
  Programming Principle~\cite[Corollary~8.3.15]{SyllaPhD}:
  \begin{eqnarray*}
    &
    & U_o^{*}(x_o) - W(x_T)
    \\
    & =
    & - \inf_{\substack{q(0) = x_o \\ q \in
    \cR_T}} \int_{0}^T L\left(q(s), \dot q(s)\right) \d{s}
    \\
    & =
    & - \inf_{\substack{q(0) = x_o \\ q \in \Lip(\mathopen]0,
    T\mathclose[;\R)}}
    \int_{0}^T L\left(q(s), \dot q(s)\right) \d{s}
    \\
    & =
    & - \inf_{\substack{q(0) = x_o, q(\tau) = \xi(\tau) \\ q \in
    \Lip(]0, \tau[;\R)}} \int_{0}^{\tau} L\left(q(s), \dot
    q(s)\right) \d{s} - \inf_{\substack{q(\tau) = \xi(\tau), q(T)
    = x_T \\ q \in \Lip(]\tau, T[;\R)}}
    \int_{\tau}^{T} L\left(q(s), \dot q(s)\right) \d{s}
    \\
    & =
    & - \int_{0}^{\tau} L\left(y(s), \dot y(s)\right) \d{s}
      - \int_{\tau}^{T} L\left(x(s), \dot x(s)\right) \d{s}
    \\
    & =
    & - \int_{0}^{T} L\left(\xi(s), \dot \xi(s)\right) \d{s} \,.
  \end{eqnarray*}
  This ensures that $\xi$ is a Lipschitz maximizer for $U_o^{*}(x_o)$,
  therefore, $\xi \in \W{2}{\infty}(\mathopen]0, T\mathclose[;\R)$
  by~\cite[Corollary~8.3.7]{SyllaPhD}. However, this contradicts the
  fact that $\xi$ is not differentiable at $t=\tau$. We conclude that
  $x$ and $y$ do not cross in $\mathopen]0, T\mathclose[$ implying
  $x(T) \leq y(T)$ and, hence, $x_T \leq y_T$.
\end{proofof}

\begin{proofof}{Theorem~\ref{th:nonemptiness}} The proof of the
  implication~\eqref{item:2} $\implies$~\eqref{item:1} is clear.

  \paragraph{Proof of~\eqref{item:1} $\implies$~\eqref{item:3}.}
  Suppose that $I_T^{HJ}(W) \neq \emptyset$ and set $w = W'$, so that
  $I_T^{CL}(w) \neq \emptyset$ by the correspondence
  between~\eqref{eq:cl} and~\eqref{eq:hj} proved
  in~\cite[Theorem~2.20]{CPS2022}.  We check that $\cG$
  in~\eqref{eq:graph} enjoys the maximal
  property~\eqref{eq:graph_maximal}. Fix $U_o \in I_T^{HJ}(W)$,
  $(x_o, x_T'), (x_o, x_T'') \in \cG$, with $x_T' < x_T''$, and
  $x_T \in \mathopen]x_T', x_T''\mathclose[$. Let $y, z \in \cR_T$, as
  defined in~\eqref{eq:hamiltonian_rays}, be solutions
  to~\eqref{eq:ode_system} connecting $(x_o, x_T')$ and
  $(x_o, x_T'')$, respectively, and let $\xi$ be the minimal backward
  generalized characteristics emanating from $(T, x_T)$, associated
  with~\eqref{eq:cl} with initial data
  $U_o'$. By~\cite[Theorem~3.2]{Dafermos1977}, $\xi$ is genuine,
  $\xi \in \cR_T$ and by Theorem~\ref{th:two_min},
  \begin{displaymath}
    W(x_T)
    = \int_0^T L\left(\xi(s), \dot \xi(s)\right) \d{s} + U_o\left(\xi(0)\right)
    \geq \int_0^T L\left(\xi(s), \dot \xi(s)\right) \d{s} + U_o^{*}\left(\xi(0)\right).
  \end{displaymath}
  Above, we used the fact that $U_o \geq U_o^{*}$, see
  Lemma~\ref{lmm:candidat_inverse_design2}.  We deduce that
  \begin{displaymath}
    U_o^{*}\left(\xi(0)\right)
    \leq
    W(x_T) - \int_0^T L\left(\xi(s), \dot \xi(s)\right) \d{s}.
  \end{displaymath}
  By definition~\eqref{eq:inverse_design} of $U_o^{*}$, we have
  equality above, and therefore $\xi$ is a point of maximum of the
  functional in~\eqref{eq:inverse_design}. We deduce that
  $(\xi(0), x_T) \in \cG$. By~\eqref{eq:graph_monotone} in
  Lemma~\ref{lmm:cs3},
  \begin{displaymath}
    x_T' < x_T \implies x_o \leq \xi(0)
    \quad \text{and} \quad x_T'' > x_T \implies x_o \geq \xi(0).
  \end{displaymath}
  We deduce that $\xi(0) = x_o$ and, therefore, $(x_o, x_T) \in \cG$.

  \paragraph{Proof of~\eqref{item:3} $\implies$~\eqref{item:2}.} We
  now show that $U_o^{*}$, as defined in~\eqref{eq:inverse_design}, is
  in $I_T^{HJ}(W)$. We first check that:
  \begin{equation}
    \label{eq:cns11}
    \forall \, x_T \in \R, \; \exists \, x_o \in \R, \quad (x_o, x_T) \in \cG.
  \end{equation}
  Note that~\eqref{eq:cns11} differs from~\eqref{item:10} in
  Lemma~\ref{lmm:cs3}, since the roles of the elements in the pair
  $(x_o, x_T)$ are reversed.

  Fix $x_T \in \R$ and introduce the subset:
  \begin{displaymath}
    E = \left\{
      x \in \R  \colon
      \exists \, y \in \mathopen]-\infty, x_T \mathclose[ \mbox{ and } (x, y) \in \cG
    \right\} \,.
  \end{displaymath}
  Fix $x \in \R$ such that $x < x_T - T \, \C{}_{H, W}$. As a
  consequence of~\eqref{item:10} in Lemma~\ref{lmm:cs3}, there exists
  $y \in \R$ such that $(x, y) \in \cG$. Now, using~\eqref{item:12} in
  Lemma~\ref{lmm:cs3}, we can write
  \begin{displaymath}
    y  = (y - x) + x
    < \modulo{x-y} + (x_T - T \, \C{}_{H, W})
    \leq x_T \,,
  \end{displaymath}
  which ensures that
  $\mathopen]-\infty, x_T - T \, \C{}_{H, W}\mathclose[ \subset E$
  and, therefore, $E$ is non-empty. Moreover, for all $x \in E$, if
  $y \in \R$ ($y < x_T$) is such that $(x, y) \in \cG$, then we have
  \begin{displaymath}
    x \leq (x - y) + y
    < \modulo{x-y} + x_T
    \leq T \, \C{}_{H, W} + x_T \,,
  \end{displaymath}
  proving that $E$ is bounded above. Hence, $\overline{x} = \sup E$ is
  finite. Likewise, the subset
  \begin{displaymath}
    F
    = \left\{ x \in \R  \colon
      \exists \, y \in \mathopen]x_T, +\infty\mathclose[
      \mbox{ and }(x, y) \in \cG
    \right\}
  \end{displaymath}
  is nonempty and bounded below. Therefore, $\underline{x} = \inf F$
  is finite. The monotonicity of $\cG$ in~\eqref{item:13} of
  Lemma~\ref{lmm:cs3} ensures that $\overline{x} \leq \underline{x}$.

  Let $(x_n)_n$ be a sequence of $E$ which converges to
  $\overline{x}$. For all $n \in \N$, there exists $y_n < x_T$ such
  that $(x_n, y_n) \in \cG$. Since $(x_n)_n$ is bounded, $(y_n)_n$ is
  bounded as well, as a consequence of~\eqref{item:12} in
  Lemma~\ref{lmm:cs3}. Up to the extraction of a subsequence, we can
  assume that $(y_n)_n$ converges to some $\overline{y} \leq
  x_T$. Since $\cG$ is closed, by~\eqref{item:10} in
  Lemma~\ref{lmm:cs3}, $(\overline{x}, \overline{y}) \in \cG$. The
  same way, there exists $\underline{y} \geq x_T$ such that
  $(\underline{x}, \underline{y}) \in \cG$. Let us conclude the proof
  by a case by case study.

  \paragraph{Case 1: $\overline{x} = \underline{x}$.} Call $x_o$ this
  common value. Since $\overline{y} \leq x_T \leq \underline{y}$, we
  have by~\eqref{eq:graph_maximal}
  \begin{displaymath}
    (x_o, \overline{y}), (x_o, \underline{y}) \in \cG
    \implies (x_o, x_T) \in \cG \,.
  \end{displaymath}

  \paragraph{Case 2: $\overline{x} < \underline{x}$.} Fix
  $x_o \in \mathopen]\overline{x},
  \underline{x}\mathclose[$. By~\eqref{item:10} in
  Lemma~\ref{lmm:cs3}, there exists $y \in \R$ such that
  $(x_o, y) \in \cG$. However, by the definition of $\overline{x}$, we
  necessarily have $y \geq x_T$. Similarly, the definition of
  $\underline{x}$ ensures that $y\leq x_T$. We proved that $y = x_T$
  and therefore, $(x_o, x_T) \in \cG$ for any
  $x_o \in \mathopen]\overline{x}, \underline{x}\mathclose[$.

  \smallskip

  Equality~\eqref{eq:cns11} rewrites as:
  \begin{equation}
    \label{eq:cns12}
    \forall \, x \in \R, \;
    \exists \, q \in \cR_T \colon
    q(T) = x
    \mbox{ and }
    W(x)
    =
    \int_{0}^T L\left(q(s), \dot q(s)\right) \d{s}
    + U_o^{*}\left(q(0)\right) \,.
  \end{equation}
  Moreover, by the definition of $U_o^{*}$, we also have
  \begin{equation}
    \label{eq:cns13}
    \forall\, q \in \mathcal{R}_T \qquad
    \int_0^T L\left(q(s), \dot q(s)\right) \d{s}
    + U_o^{*}\left(q(0)\right)
    \geq
    W\left(q(T)\right) \,.
  \end{equation}
  Together~\eqref{eq:cns12} and~\eqref{eq:cns13} imply that
  \begin{displaymath}
    \begin{aligned}
      \forall x \in \R, \quad W(x) & = \inf_{\substack{q(T) = x \\ q
          \in \cR_T}} \int_{0}^T L\left(q(s), \dot q(s)\right) \d{s} +
      U_o^{*}\left(q(0)\right)
      \\
      & = \inf_{\substack{q(T) = x \\ q \in \Lip(\mathopen]0,
          T\mathclose[;\R)}} \int_{0}^T L\left(q(s), \dot q(s)\right)
      \d{s} + U_o^{*}\left(q(0)\right) \,,
    \end{aligned}
  \end{displaymath}
  by~\cite[Corollary~8.3.15]{SyllaPhD}. This last equality means that
  the viscosity solution $U$ to~\eqref{eq:hj} associated with initial
  datum $U_o^{*}$ verifies $U(T) = W$, using the classical
  correspondence viscosity solution/calculus of variations,
  see~\cite[Theorem~8.3.12]{SyllaPhD}. We proved that
  $U_o^{*} \in I_T^{HJ}(W)$.

  \paragraph{Proof of~\eqref{item:1} $\implies$ $\pi_{W'}$ is well
    defined and nondecreasing.} Suppose that
  $I_T^{HJ}(W) \neq \emptyset$ and set $w = W'$, so that
  $I_T^{CL}(w) \neq \emptyset$ by~\cite[Theorem~2.20]{CPS2022}. In the
  light of both Lemma~\ref{lmm:regularity_profile_cl} and
  Lemma~\ref{lem:HS}, $\pi_w$ is well-defined by~\eqref{eq:4}.

  Fix $x, y \in \R$ with $x < y$. Since $I_T^{CL}(w) \neq \emptyset$,
  $\pi_w$ assigns to $x$, respectively $y$, the value at time $t=0$ of
  the minimal backward generalized characteristics emanating from
  $(T, x)$, respectively from $(T, y)$, which we denote by $\xi_x$,
  respectively $\xi_y$.  By~\cite[Theorem~3.2]{Dafermos1977}, $\xi_x$
  and $\xi_y$ are genuine, hence they do not intersect in
  $\mathopen]0, T\mathclose[$,
  see~\cite[Corollary~3.2]{Dafermos1977}. This implies in particular
  that $\xi_x(0) \leq \xi_y(0)$, proving that $\pi_w$ is
  nondecreasing.
\end{proofof}

\subsection{Proof of Theorem~\ref{th:belonging}}
\label{sec:proofs-relat-sect}

\begin{proofof}{Theorem~\ref{th:belonging}}
  We prove the two implications separately.

  \paragraph{Claim: If $U_o \in I^{HJ}_T (W)$, then \emph{(i)}
    and~\emph{(ii)} hold.}

  Point \textit{(i)} comes from
  Lemma~\ref{lmm:candidat_inverse_design2}. Let us prove that
  \textit{(ii)} holds. Fix $x_o \in \pi_{W'}(\R)$. By definition,
  there exists an $x \in \R$ such that $x_o = \pi_{W'}(x)$.  This
  means that $x_o$ is the value at time $t=0$ of the minimal backward
  characteristics $\xi$, see~\cite[Definition~3.1, Theorems~3.2
  and~3.3]{Dafermos1977} emanating from $(T, x)$. Since
  $U_o \in I_T^{HJ}(W)$, Theorem~\ref{th:two_min} ensures that
  \begin{displaymath}
    W(x) = \int_0^T L\left(\xi(s), \dot \xi(s)\right) \d{s} + U_o(x_o) \,.
  \end{displaymath}
  On the other hand, by~\eqref{eq:inverse_design},
  \begin{displaymath}
    U_o^{*}(x_o)
    =
    \sup_{\substack{y(0) = x_o \\ y \in \cR_T}}
    \left( W\left(y(T)\right) {-} \int_{0}^T \! L\left(y(s), \dot y(s)\right) \right) \d{s}
    \geq
    W(x) - \int_0^T \! L\left(\xi(s), \dot \xi(s)\right) \d{s}
    =
    U_o(x_o) .
  \end{displaymath}
  So, $U_o = U_o^{*}$ on $\pi_{W'}(\R)$. These two functions are
  continuous, hence they coincide on $\overline{\pi_{W'}(\R)}$.

  \paragraph{Claim: If \emph{(i)} and~\emph{(ii)} hold, then
    $U_o \in I^{HJ}_T (W)$.}

  Fix $x \in \R$.  Recall that $I_T^{HJ}(W) \neq \emptyset$ which, by
  Theorem~\ref{th:nonemptiness}, ensures that
  $U_o^{*} \in I_T^{HJ}(W)$. This, together with the inequality
  $U_o \geq U_o^{*}$, immediately implies:
  \begin{eqnarray}
    \nonumber
    W(x)
    & =
    & \inf_{\substack{y(T) = x \\y \in \cR_T}}
    \left(
    \int_{0}^T L\left(y(s), \dot y(s)\right) \d{s} + U_o^*\left(y(0)\right)
    \right)
    \\
    \label{eq:cns21}
    & \leq
    & \inf_{\substack{y(T) = x \\y \in \cR_T}}
    \left( \int_{0}^T L\left(y(s), \dot y(s)\right) \d{s} + U_o\left(y(0)\right) \right).
  \end{eqnarray}
  Denote by $\xi$ the minimal backward characteristics emanating from
  $(T, x)$. Using both the facts that $U_o^{*} \in I_T^{HJ}(W)$ and
  that, by Theorem~\ref{th:two_min}, $\xi$ is a minimizer, we have:
  \begin{displaymath}
    W(x) = \int_0^T L\left(\xi(s), \dot \xi(s)\right) \d{s} + U_o^{*}\left(\xi(0)\right).
  \end{displaymath}
  Clearly, $\xi(0) \in \overline{\pi_{W'}(\R)}$ and therefore,
  by~\emph{(ii)}, we can replace $U_o^{*}\left(\xi(0)\right)$ by
  $U_o\left(\xi(0)\right)$ in the last equality. This ensures that we
  have equality in~\eqref{eq:cns21}, which means that
  $U_o \in I_T^{HJ}(W)$.
\end{proofof}

\subsection{Proof of Theorem~\ref{th:counter_example}}
\label{subsec:proof-theor-refth:c}

In all proofs in this section, the reader might want to keep
Figure~\ref{fig:nice} in mind for a helpful geometrical visualization.

Long but straightforward computations show that $H$, as defined
in~\eqref{eq:11}, satisfies~\eqref{eq:smoothness},
\eqref{eq:space_dpcy} with $X = 1$, and~\eqref{eq:convexity}, see
Figure~\ref{fig:g}. With this flux, the conservation law~\eqref{eq:cl}
is also the inviscid Burger equation with source term $-g'$, see
Figure~\ref{fig:g}. We fix the initial datum
\begin{equation}
  \label{eq:12}
  u_o(x) \coloneqq
  \left\{
    \begin{array}{cl}
      -2 & \text{if } x < 0\,,
      \\
      2 & \text{if } x > 0\,,
    \end{array}
  \right.
\end{equation}
which would evolve into a rarefaction in the homogeneous case. The
proof of Theorem~\ref{th:counter_example} is based on the Cauchy
problem for~\eqref{eq:ode_system} which, in this case, reads
\begin{equation}
  \label{eq:ode_system_ce}
  \left\{
    \begin{array}{rcl}
      \dot q
      & =
      & p
      \\
      \dot p
      & =
      & -g'(q)
      \\
      q(0)
      & =
      & q_o
      \\
      p(0)
      & =
      & p_o
    \end{array}
  \right.
  \qquad \mbox{with } g \mbox{ as in~\eqref{eq:11},}
\end{equation}
and to which Lemma~\ref{lem:HS} applies.  The first equation
in~\eqref{eq:ode_system_ce} will be tacitly used throughout this
section. By the Hamiltonian nature of~\eqref{eq:ode_system_ce}, $H$ is
conserved along solutions, so that
\begin{equation}
  \label{eq:conservation_ce}
  \forall t \in \R, \quad \frac{p(t)^2}{2} + g\left(q(t)\right) = \frac{p_o^2}{2} + g(q_o) \,.
\end{equation}

\begin{lemma}
  \label{lmm:ce1}
  Let $H$ be as in~\eqref{eq:11} and $u_o$ be as in~\eqref{eq:12}.
  Fix $q_o \geq 0$. Denote by $(q, p)$ the solution
  to~\eqref{eq:ode_system_ce} with initial datum
  $\left(q_o, u_o(q_o+)\right) = (q_o, 2)$. Then, $q$ is increasing on
  $[0, +\infty\mathclose[$ and $q(t) \limit{t}{+\infty} +\infty$.
\end{lemma}

\begin{proofof}{Lemma~\ref{lmm:ce1}}
  Note that $p_o >0$.  By~\eqref{eq:conservation_ce}, for all
  $t \in \R$
  \begin{displaymath}
    p (t)^2
    =
    \underbrace{{p (0)}^2}_{=4} + 2 \underbrace{g (q_o)}_{\geq 0} - 2 \underbrace{g\left(q (t)\right)}_{\leq 1}
    \geq
    2 \,.
  \end{displaymath}
  Thus, for $t \in \R$, $p (t) {\geq}
  \sqrt2$. By~\eqref{eq:ode_system}, $q$ is strictly increasing and
  $q (t) {\geq} \sqrt2 \, t$ for $t {\in} [0, +\infty\mathclose[$.
\end{proofof}

\begin{figure}[!ht]
  \begin{center}
    \includegraphics[scale=0.5]{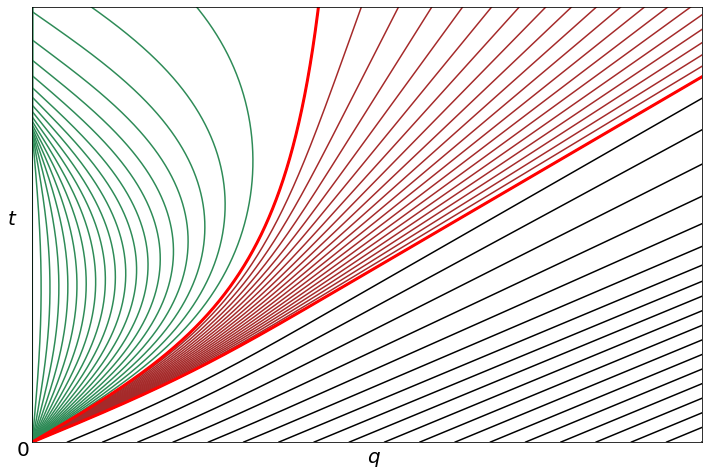}
  \end{center}
  \caption{On the horizontal axis, the $q$ component of solutions
    to~\eqref{eq:ode_system_ce}, while $t$ is on the vertical
    axis. Brown curves are those considered in~\eqref{item:17} of
    Lemma~\ref{lmm:ce3}; green curves refer to
    Lemma~\ref{lmm:reallyLast}. The 2 red thicker curves depict
    solutions corresponding to the initial data $(0, \sqrt{2})$ and
    $(0, 2)$. The black curves are those considered in
    Lemma~\ref{lmm:ce1} and in Lemma~\ref{lmm:ce2}.}
  \label{fig:fc2}
\end{figure}
Refer to the lines on the right in Figure~\ref{fig:fc2} for an
illustration of the different behaviors of $q$ described in
Lemma~\ref{lmm:ce1} and Lemma~\ref{lmm:ce2}.

\begin{lemma}
  \label{lmm:ce2}
  Let $H$ be as in~\eqref{eq:11} and $u_o$ be as in~\eqref{eq:12}.
  Fix $0 \leq q_o < \widetilde{q_o}$ and denote by $(q, p)$,
  respectively $(\widetilde{q}, \widetilde{p})$, the global solution
  to~\eqref{eq:ode_system} with initial datum
  $\left(q_o, u_o(q_o+)\right) = (q_o,2)$, respectively
  $\left(\widetilde{q_o}, u_o (\widetilde{q_o})\right) =
  (\widetilde{q_o}, 2)$. Then, $q(t) < \widetilde{q}(t)$, for all
  $t \geq 0$.
\end{lemma}

\begin{proofof}{Lemma~\ref{lmm:ce2}}
  Set $p_o = u_o(q_o)$ and $\widetilde{p_o} =
  u_o(\widetilde{q_o})$. We proceed by contradiction. Let $\tau > 0$,
  be the smallest time where $q (\tau) = \tilde q (\tau)$. Since
  $q_o < \tilde q_o$, we have that $p (\tau) \geq \tilde p (\tau)$. By
  Lemma~\ref{lmm:ce1}, $p(\tau) \geq \widetilde{p}(\tau) \geq 0$.
  Then,
  \begin{displaymath}
    \begin{aligned}
      \left.
        \begin{array}{@{}r@{\,}c@{\,}l@{}}
          q_o
          & \in
          & [0,\widetilde{q_o}\mathclose[
          \\
          p_o
          & =
          & \tilde p_o
        \end{array}
      \right\} & \implies \frac{p_o^2}{2} + g(q_o) \leq
      \frac{\widetilde{p_o}^2}{2} + g(\widetilde{q_o}) &
      [\mbox{By~\eqref{eq:11} and~\eqref{eq:12}}]
      \\
      & \implies \frac{p(\tau)^2}{2} + g\left(q(\tau)\right) \leq
      \frac{\widetilde{p}(\tau)^2}{2} +
      g\left(\widetilde{q}(\tau)\right) &
      [\mbox{By~\eqref{eq:conservation_ce}}]
      \\
      & \implies p(\tau) \leq \widetilde{p}(\tau).  & [\mbox{By }
      q(\tau) = \widetilde{q}(\tau),\; p(\tau) \geq
      \widetilde{p}(\tau) > 0]
    \end{aligned}
  \end{displaymath}
  We then deduce that $p(\tau) = \widetilde{p}(\tau)$, which
  contradicts the uniqueness proved by Cauchy Lipschitz Theorem.
\end{proofof}

\begin{lemma}
  \label{lmm:ce3}
  Let $H$ be as in~\eqref{eq:11} and $u_o$ be as in~\eqref{eq:12}. Fix
  $p_o \in [\sqrt2, 2\mathclose[$. Denote by $(q, p)$ the global
  solution to~\eqref{eq:ode_system_ce} with initial datum $(0,
  p_o)$. (Refer to Figure~\ref{fig:fc2}.)
  \begin{enumerate}[(i)]

  \item\label{item:17} If $p_o \in \mathopen]\sqrt{2}, 2\mathclose[$,
    then $q$ is increasing on $[0, +\infty\mathclose[$ and
    $q(t) \limit{t}{+\infty} +\infty$.

  \item\label{item:19} If $p_o = \sqrt{2}$, then $q$ is increasing on
    $[0, +\infty\mathclose[$, $q(t) \limit{t}{+\infty} 1$ and $q$ is
    concave.

  \end{enumerate}
\end{lemma}

\noindent Refer to the middle curves in Figure~\ref{fig:fc2} for an
illustration of the different behaviors of $q$ described in
Lemma~\ref{lmm:ce3}.

\begin{proofof}{Lemma~\ref{lmm:ce3}}
  The proof of~\eqref{item:17} is identical to that of
  Lemma~\ref{lmm:ce1}, so we omit it.

  Concerning~\eqref{item:19}, $p$ is positive on
  $\mathopen]0, +\infty\mathclose[$. Indeed, assume by contradiction
  that there exists a minimal $\tau > 0$ such that $p(\tau) =
  0$. By~\eqref{eq:conservation_ce}, we deduce that $q(\tau) =
  1$. However, $(q_s,p_s) \colon t \mapsto (1,0)$ is the unique global
  solution to~\eqref{eq:ode_system_ce} with datum $(1, 0)$. Hence, $p$
  is positive on $\mathopen]0, +\infty\mathclose[$.

  Thus, by~\eqref{eq:ode_system_ce}, $q$ is increasing on
  $[0, +\infty\mathclose[$ and positive on $]0,
  +\infty\mathclose[$. Once again, \eqref{eq:conservation_ce} and the
  presence of the stationary solution $(q_s,p_s)$ ensure that for all
  $t \in [0, +\infty\mathclose[$, $q(t) < 1$. Therefore, as
  $t \to +\infty$, $q$ admits a finite limit, say $q_\infty$, which is
  not greater than $1$.

  Moreover, the positivity of $q$ ensures,
  by~\eqref{eq:ode_system_ce}, that $p$ is nonincreasing on
  $[0, +\infty\mathclose[$, so that by~\eqref{eq:ode_system_ce}, $q$
  is concave. Since $p$ is also bounded, by~\eqref{eq:3} in
  Lemma~\ref{lem:HS}, $p$ admits a finite limit as $t \to
  +\infty$. Hence, $\dot q$ has a finite limit as $t \to +\infty$ and,
  since we already showed that $q$ converges to $q_\infty$ as
  $t \to +\infty$, then $\dot q \to 0$ and therefore $p\to 0$ as
  $t \to +\infty$. By~\eqref{eq:conservation_ce}, we get
  $g (q_\infty) = 1$, and hence $q_\infty \geq 1$, therefore
  $q_\infty = 1$, completing the proof of~\eqref{item:19}.
\end{proofof}

\begin{lemma}
  \label{lmm:reallyLast}
  Let $H$ be as in~\eqref{eq:11} and $u_o$ be as in~\eqref{eq:12}. Let
  $p_o \in \mathopen]0, \sqrt{2}\mathclose[$. Denote by $(q, p)$ the
  global solution to~\eqref{eq:ode_system_ce} with initial data
  $(0, p_o)$. Then, $q$ is periodic. Introduce the map
  \begin{equation}
    \label{eq:15}
    \fonction{\mathcal{T}}{\mathopen]0, \sqrt2\mathclose[}{\mathclose]0,+\infty\mathclose[}{p_o}{
      \text{the smallest period of $q$}}
  \end{equation}
  \begin{enumerate}[(i)]
  \item\label{item:28} $q$ is concave on $[0, \mathcal{T} (p_o)/2]$.
  \item \label{item:29} For all
    $t \in \left[0, \mathcal{T} (p_o)/2\right]$,
    $q (t) = q \left(\mathcal{T} (p_o)/2-t\right)$.
  \item \label{item:30} For all
    $t \in \left[0, \mathcal{T} (p_o)\right]$,
    $q (t) = - q \left(\mathcal{T} (p_o)-t\right)$.
  \item \label{item:31} $q$ admits its maximum at
    $\mathcal{T} (p_o)/4$.
  \end{enumerate}
\end{lemma}

\begin{proofof}{Lemma~\ref{lmm:reallyLast}}
  Note first that by~\eqref{eq:conservation_ce}, $q$ is bounded, since
  for all $t \in \R$,
  $g\left(q(t)\right) = \frac{p_o^2}{2} -
  \frac{\left(p(t)\right)^2}{2} \leq \frac{p_o^2}{2} < 1$ and hence
  $\modulo{q(t)} < 1$ for all $t \in \R$.

  Assume now, by contradiction, that $q$ does not vanish on
  $\mathopen]0,+\infty\mathclose[$. Since $q (0) = 0$ and
  $\dot q(0) > 0$, we have that for all
  $t \in \mathopen]0,+\infty\mathclose[$, $q(t) > 0$. Therefore, $p$
  is decreasing on $[0, +\infty\mathclose[$
  by~\eqref{eq:ode_system_ce}, bounded by~\eqref{eq:conservation_ce}
  and thus admits a finite limit as $t \to +\infty$.

  Thus, $q$ is bounded and its derivative $\dot q =p$ has a finite
  limit as $t \to +\infty$, hence $\lim_{t\to +\infty} \dot q (t) = 0$
  and also $\lim_{t\to+\infty} p (t) = 0$. Therefore, by monotonicity,
  $p$ is nonnegative on all $[0, +\infty\mathclose[$. Consequently,
  $q$ is nondecreasing and bounded, therefore it admits a finite limit
  $q_\infty \geq 0$ as $t \to +\infty$. On the one hand, by taking the
  limit as $t \to +\infty$ in~\eqref{eq:conservation_ce}, we get:
  \begin{displaymath}
    g(q_\infty) = \frac{p_o^2}{2} \in \mathopen]0, 1\mathclose[
    \implies q_\infty \in \mathopen]0,1\mathclose[
    \implies g'(q_\infty) \neq 0.
  \end{displaymath}
  On the other hand,
  $\dot p(t) = - g'\left(q(t)\right) \to - g'(q_\infty)$. Moreover,
  $p$ has a finite limit as $t \to +\infty$, $\dot p (t) \to 0$. This
  provides the needed contradiction.

  Thus we proved that there exists $\tau > 0$ such that $q(\tau) =
  0$. As a consequence, the number
  \begin{equation}
    \label{eq:14}
    \tau_*
    \coloneqq
    \sup \left\{t
      \in \mathopen]0, +\infty\mathclose[ \colon
      \forall \, s \in \mathopen]0, t \mathclose[ ,\;  q(s) > 0
    \right\}
  \end{equation}
  is well-defined and satisfies ${q ( \tau_*) = 0}$. Note that for
  $t \in \mathopen]0,\tau_*\mathclose[$, $q$ is positive, $p$ is
  decreasing and hence $q$ is concave on $[0, \tau_*]$,
  proving~\eqref{item:28} on $[0, \tau_*]$. Furthermore,
  $p(\tau_*) < p_o$. Apply~\eqref{eq:conservation_ce} at time
  $t = \tau_*$ to obtain $p(\tau_*)^2 = p_o^2$, which implies
  $p(\tau_*) = -p_o$.

  Now we verify that $q$ is $2\tau_*$--periodic. To this aim,
  introduce $\xi(t) = -q(t + \tau_*)$ and $\nu(t) = -p(t + \tau_*)$.
  Thanks to $g$ being even, it is straightforward to check that both
  $(q, p)$ and $(\xi, \nu)$ solve the same Cauchy
  problem~\eqref{eq:ode_system_ce}. Consequently,
  \begin{equation}
    \label{eq:13}
    \forall \, t \in \R \qquad
    q(t + 2\tau_*) = -q(t + \tau_*) = q(t)\,.
  \end{equation}
  Hence, $q$ is $2\tau_*$--periodic. By~\eqref{eq:14}, $2\tau_*$ is
  the minimal period, completing the proof of~\eqref{item:28}.

  Finally, define $\widehat{q} (t) = q (\tau_*-t)$ and
  $\widehat{p} (t) = -p (\tau_*-t)$.  Note that $(q,p)$ and
  $(\widehat{q},\widehat{p})$ both solve~\eqref{eq:ode_system_ce} with
  datum $(0, p_o)$, since $p(\tau_*) = -p_o$. Hence, for all
  $t \in [0, \tau_*]$, $q (t) = q (\tau_*-t)$,
  proving~\eqref{item:29}. Combined with the concavity of $q$, this
  ensures that $\max_{t \in [0, \tau_*]} q (t) = q (\tau_*/2)$,
  proving~\eqref{item:31}. Finally, \eqref{eq:13} and~\eqref{item:29}
  imply~\eqref{item:30}.
\end{proofof}

\begin{lemma}
  \label{lem:Last}
  Let $H$ be as in~\eqref{eq:11} and $u_o$ be as
  in~\eqref{eq:12}. Call $\mathcal{T}$ the map defined
  in~\eqref{eq:15}. Then:
  \begin{enumerate}[(i)]
  \item \label{item:24} $\mathcal{T}$ is continuous.
  \item \label{item:25} $\mathcal{T}$ strictly increasing.
  \item \label{item:26}
    $\inf_{p_o \in \mathopen]0, \sqrt{2}\mathclose[}\mathcal{T} (p_o)
    = \left.\pi\middle/\sqrt2 \right.$.
  \item \label{item:27}
    $\lim_{p_o \to \sqrt2} \mathcal{T} (p_o) = +\infty$.
  \end{enumerate}
\end{lemma}

\begin{proofof}{Lemma~\ref{lem:Last}}
  \paragraph{Proof of~\eqref{item:24}.}
  Fix $p_o \in \mathopen]0, \sqrt2\mathclose[$ and let $q_o = 0$. Let
  $(q,p)$ be the solution to~\eqref{eq:ode_system_ce}. Then, by
  Lemma~\eqref{lmm:reallyLast}, we get $p (t) = \dot q (t) >0$ for all
  $t \in \mathopen]0, \mathcal{T} (p_o)/4\mathclose[$. Hence,
  using~\eqref{eq:conservation_ce} we get
  \begin{displaymath}
    \forall \, t \in \mathopen]0, \mathcal{T} (p_o) /4\mathclose[ \qquad
    p (t) = \sqrt{p_o^2 - 2 g\left(q (t)\right)}
  \end{displaymath}
  and by~\eqref{eq:ode_system_ce}, we have
  \begin{equation}
    \label{eq:23}
    \frac{\mathcal{T} (p_o)}{4} = \int_0^{\mathcal{T} (p_o)/4} \frac{\dot q(t)}{
      \sqrt{p_o^2 - 2 g\left(q(t)\right)}} \d{t}  \,.
  \end{equation}
  Note that the integrand in the right hand side above is singular
  when $t = \mathcal{T} (p_o) /4$, but it is positive for all $t$.
  Use the change of variable $x = q (t)$ to get
  \begin{equation}
    \label{eq:28}
    \frac{\mathcal{T} (p_o)}{4}
    =
    \int_0^{g^{-1}(p_o^2 / 2)}
    \frac{1}{\sqrt{p_o^2 - 2 g(x)}} \d{x}
  \end{equation}
  where $g^{-1}$ is the inverse of the $\C1$ diffeomorphism
  $g_{\strut|\mathopen]0,1\mathopen[} \colon \mathopen]0,1\mathopen[
  \to \mathopen]0,1\mathopen[$.

  Define $\mathcal{A}\colon \mathopen]0,1\mathopen[ \to \R_+$ by
  \begin{equation}
    \label{eq:26}
    \mathcal{A} (r)
    \coloneqq
    \int_0^1 \dfrac{r}{\sqrt{g (r) - g (\theta \, r)}} \d\theta \,,
  \end{equation}
  so that the change of variable $x = \theta \, r$ with
  $r = g^{-1} (p_o^2/2)$ in~\eqref{eq:28} leads to
  \begin{equation}
    \label{eq:27}
    \mathcal{T} (p_o)
    =
    2 \sqrt 2 \; \mathcal{A}\left(g^{-1} (p_o^2/2)\right) \,.
  \end{equation}
  The continuity of $\mathcal{A}$ is proved in~Lemma~\ref{lem:NoHope}
  in Appendix~\ref{sec:appendix}, completing the proof
  of~\eqref{item:24}.

  \paragraph{Proof of~\eqref{item:25}.}
  By~\cite[Theorem~A]{Chicone1987}, \eqref{eq:27}
  and~\eqref{eq:conservation_ce}, the condition
  \begin{displaymath}
    \forall x \in \mathopen]0, 1\mathclose[, \quad \frac{\d{}^2~}{\d{x}^2}\left(\frac{g(x)}{g'(x)^2}\right) \geq 0
  \end{displaymath}
  ensures that
  $E \mapsto 2\sqrt 2\, \mathcal{A}\left(g^{-1} (E)\right)$ is
  increasing. Hence, by~\eqref{eq:27}, also $\mathcal{T}$ is
  increasing on $\mathopen]0, \sqrt2\mathclose[$. By~\eqref{eq:11},
  for $x \in \mathopen]0,1\mathclose[$ we have
  \begin{displaymath}
    \frac{\d{}^2~}{\d{x}^2}\left(\frac{g(x)}{g'(x)^2}\right)
    = -
    \frac{3}{32} \frac{7x^8 -32 x^6 + 59x^4 -56 x^2 -6}{(x^2 - 1)^8} \,.
  \end{displaymath}
  We leave to Lemma~\ref{lem:Sturm} in Appendix~\ref{sec:appendix} the
  proof that
  $\frac{\d{}^2~}{\d{x}^2}\left(\frac{g(x)}{g'(x)^2}\right) \geq 0$
  for all $x \in \mathopen]0, 1\mathclose[$ by means of Sturm Theorem,
  see~\cite{Sturm1829}.

  \paragraph{Proof of~\eqref{item:26}.}
  To prove the lower bound on $\mathcal{T}$, introduce for any
  $p_o \in \mathopen]0, \sqrt2\mathclose[$, $\tilde q>0$ so that
  $g(\tilde q) = p_o^2/2$:
  \begin{flalign*}
    \inf_{p_o \in ]0, \sqrt2]} \mathcal{T} (p_o) %
    & = \lim_{p_o \to 0+} \mathcal{T} (p_o) %
    & [\mbox{By the monotonicity of $\mathcal{T}$}]
    \\
    & = \lim_{E \to 0}2 \sqrt 2 \; \mathcal{A}\left(g^{-1} (E)\right)
    & [\mbox{By~\eqref{eq:27}}]
    \\
    & = \lim_{r\to 0} 2 \sqrt 2 \; \mathcal{A}(r) &
    [\mbox{By~\eqref{eq:11}}]
    \\
    & = 2\sqrt2 \lim_{r \to 0} \int_0^1 \frac{r}{\sqrt{g (r) - g
        (\theta \, r)}} \d\theta %
    & [\mbox{By~\eqref{eq:26}}]
    \\
    & = 4 \int_0^1 \frac{\d\theta}{\sqrt{g'' (0) \, (1-\theta^2)}}
    \\
    & = \frac{\pi}{\sqrt2} \,.
  \end{flalign*}
  This completes the proof of~\eqref{item:26}.

  \paragraph{Proof of~\eqref{item:27}.}
  Similar computations, using now Fatou's Lemma, lead to
  \begin{displaymath}
    \lim_{p_o \to \sqrt2} \mathcal{T} (p_o) %
    = 2\sqrt2 \lim_{r \to 1} \int_0^1 \frac{r}{\sqrt{g (r) - g
        (\theta \, r)}} \d\theta
    \geq 2\sqrt 2 \int_0^1 \frac{1}{(1-\theta^2)^2} \d{\theta}
    = +\infty \,,
  \end{displaymath}
  completing the proof of~\eqref{item:27} and of Lemma~\ref{lmm:ce3}.
\end{proofof}

\begin{lemma}
  \label{cor:ce4}
  Let $H$ be as in~\eqref{eq:11} and $u_o$ be as in~\eqref{eq:12}.
  Fix $0 < p_o < \widetilde{p_o} < 2$ and denote by $(q, p)$,
  respectively $(\widetilde{q}, \widetilde{p})$, the global solution
  to~\eqref{eq:ode_system_ce} with initial datum $(0, p_o)$,
  respectively $(0, \widetilde{p_o})$. Then,
  \begin{displaymath}
    \begin{array}{c@{\quad}r@{\,}c@{\,}lcr@{\,}c@{\,}l@{\qquad}r@{\,}c@{\,}l}
      (i)
      & p_o
      & \in
      & \mathopen]0, \sqrt2\mathclose[
      & \implies
      & \forall\, t
      & \in
      & \mathopen]0, \mathcal{T} (p_o)/2]
      & q (t)
      & <
      & \widetilde{q} (t)\,;
      \\
      (ii)
      & p_o
      & \in
      & [\sqrt2,2\mathclose[
      & \implies
      & \forall\, t
      & \in
      & \mathopen]0, +\infty\mathclose[
      & q (t)& <
      & \widetilde{q} (t) \,.
    \end{array}
  \end{displaymath}
\end{lemma}

\noindent Refer to the middle and left curves in Figure~\ref{fig:fc2}
for an illustration of the different behaviors of $q$ described in
Lemma~\ref{cor:ce4}.

\begin{proofof}{Lemma~\ref{cor:ce4}}
  We split the proof in several steps.

  \paragraph{Claim~1: Let $T>0$ be such that $p (t)>0$ and
    $\widetilde{p} (t) >0$ for all $t \in [0,T]$. Then, for all
    $t \in [0,T]$, $q (t) < \widetilde{q} (t)$.}
  By contradiction, since $p_o < \widetilde{p}_o$, there exists
  $s \in \mathopen]0,T]$ such that $q (t) < \widetilde{q} (t)$ for
  $t \in \mathopen]0, s\mathclose[$, $q (s) = \widetilde{q} (s)$ and
  thus $p (s) \geq \widetilde{p} (s)$. Then,
  by~\eqref{eq:conservation_ce},
  \begin{flalign*}
    0 < p_o < \widetilde{p}_o %
    & \implies H (0,p_o) < H (0,\widetilde{p}_o) %
    & [\mbox{By~\eqref{eq:11}}]
    \\
    & \implies H \left(q(s),p(s)\right) < H\left(\widetilde{q} (s),
      \widetilde{p} (s)\right)%
    & [\mbox{By~\eqref{eq:conservation_ce}}]
    \\
    & \implies p (s) < \widetilde{p} (s) & [\mbox{Since } q (s) =
    \widetilde{q} (s) \mbox{ and } p (s)>0,\;\widetilde{p} (s)>0]
  \end{flalign*}
  which yields a contradiction, proving Claim~1.

  \paragraph{Claim~2: \emph{(i)} holds for
    $\widetilde{p}_o \in \mathopen]0, \sqrt2 \mathclose[$.}
  By Claim~1, for all
  $t \in \mathopen]0, \mathcal{T} (p_o)/4\mathclose[$,
  $q (t) < \widetilde{q} (t)$. Indeed, by~\eqref{eq:ode_system_ce}
  together with Lemma~\ref{lmm:reallyLast} and Lemma~\ref{lem:Last},
  both $p$ and $\widetilde{p}$ are positive on
  $\mathopen]0, \mathcal{T} (p_o)/4\mathclose[$ and Claim~1 applies.
  Hence, by the symmetry in~\eqref{item:29} of
  Lemma~\ref{lmm:reallyLast}, we have
  \begin{displaymath}
    \forall\ t \in
    \left]
      \frac{\mathcal{T} (\widetilde{p}_o)}{2} - \frac{\mathcal{T} (p_o)}{4}, \frac{\mathcal{T} (\widetilde{p}_o)}{2} \right[ \qquad
    q \left(\frac{\mathcal{T} (\widetilde{p}_o)}{2} -t \right)
    < \widetilde{q} (t) \,.
  \end{displaymath}
  Introduce the concave function
  \begin{equation}
    \label{eq:16}
    \fonction{\eta}{\left[0,\frac{\mathcal{T} (\widetilde{p}_o)}{2}\right]}{\R}{t}{%
      \left\{
        \begin{array}{l@{\qquad}r@{\,}c@{\,}l}
          q (t)
          & t
          & \in
          & \left[0,\frac{\mathcal{T} (p_o)}{4}\right]
          \\[6pt]
          q \left(\frac{\mathcal{T} (p_o)}{4}\right)
          & t
          & \in
          & \left]\frac{\mathcal{T} (p_o)}{4}, \frac{\mathcal{T} (\widetilde{p}_o)}{2} - \frac{\mathcal{T} (p_o)}{4}\right[
          \\[6pt]
          q \left(\frac{\mathcal{T} (\widetilde{p}_o)}{2} -t \right)
          & t
          & \in
          & \left[
            \frac{\mathcal{T} (\widetilde{p}_o)}{2} - \frac{\mathcal{T} (p_o)}{4}, \frac{\mathcal{T} (\widetilde{p}_o)}{2} \right]
        \end{array}
      \right.
    }
  \end{equation}
  and note that, see Figure~\ref{fig:concavity},
  \begin{figure}[!ht]
    \begin{center}
      \includegraphics[scale=0.5]{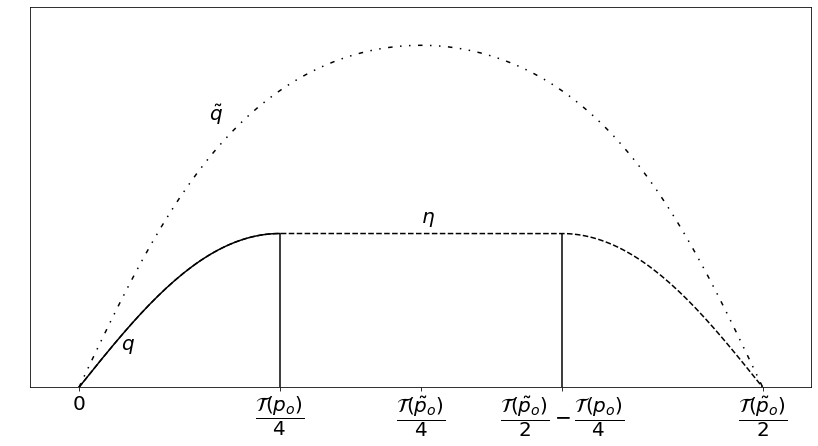}
    \end{center}
    \caption{Curves used in Claim~2 in the proof of
      Lemma~\ref{cor:ce4}. The dashed curve is the graph of $\eta$
      in~\eqref{eq:16}. The continuous curves are the graphs of $q$
      restricted to $[0, \mathcal{T} (p_o)/2]$ and of its
      translate. The dashed--dotted curved is the graph of
      $\widetilde{q}$.}
    \label{fig:concavity}
  \end{figure}
  \begin{displaymath}
    \begin{array}{r@{\,}c@{\,}l@{\qquad}r@{\,}c@{\,}l@{\qquad}l}
      \forall\, t
      & \in
      & \left[0, \frac{\mathcal{T} (p_o)}{4}\right]
      & q (t)
      & \leq
      & \eta (t)
      & [\mbox{By~\eqref{eq:16}}]
      \\[6pt]
      \forall\, t
      & \in
      & \left[0, \frac{\mathcal{T} (\widetilde{p}_o)}{2}\right]
      & \eta (t)
      & \leq
      & \widetilde{q} (t)
      & [\mbox{By concavity of $\eta$ and $\widetilde{q}$}]
      \\[6pt]
      \forall\, t
      & \in
      & \left[\frac{\mathcal{T} (p_o)}{4}, \frac{\mathcal{T} (p_o)}{2} \right]
      & q (t)
      & \leq
      & \eta (t)
      & [\mbox{By symmetry}]
    \end{array}
  \end{displaymath}
  completing the proof of Claim~2.

  \paragraph{Claim~3: \emph{(i) holds for
      $\widetilde{p}_o \in [\sqrt2, +\infty\mathclose[$.}} By Claim~1,
  for all $t \in \mathopen]0, \mathcal{T} (p_o)/4\mathclose]$,
  $q (t) < \widetilde{q} (t)$.

  By~\eqref{item:31} in Lemma~\ref{lmm:reallyLast}, for all
  $t \in \mathopen]0, \mathcal{T} (p_o)/2\mathclose[$,
  $q (t) \leq q (\mathcal{T} (p_o)/4)$, while by~\eqref{item:17}
  and~\eqref{item:19} in Lemma~\ref{lmm:ce3}, for
  $\widetilde{p}_o \in [\sqrt2, +\infty\mathclose[$ and for all
  $t \in [\mathcal{T} (p_o)/4, +\infty\mathclose[$,
  $\widetilde{q} (t) > \widetilde{q} (\mathcal{T} (p_o)/4)$.  All this
  ensures that
  \begin{displaymath}
    \forall \, t \in
    \left[\frac{\mathcal{T} (p_o)}{4}, \frac{\mathcal{T} (p_o)}{2}\right]
    \qquad
    q (t)
    \leq q \left(\frac{\mathcal{T} (p_o)}{4}\right)
    < \widetilde{q} \left(\frac{\mathcal{T} (p_o)}{4}\right)
    \leq \widetilde{q} (t) \,,
  \end{displaymath}
  completing the proof of Claim~3.

  \paragraph{Claim~4: Proof of~\emph{(ii)}.} If
  $p_o, \widetilde{p_o} \in [\sqrt{2}, 2[$, then by~\eqref{item:17}
  and~\eqref{item:19} in Lemma~\ref{lmm:ce3}, $q$ and $\widetilde{q}$
  are increasing, so that by~\eqref{eq:ode_system_ce} $p$ and
  $\widetilde{p}$ are positive. So, Claim~1 applies, completing the
  proof.

\end{proofof}


We use below the flow $\flow$ introduced in~\eqref{eq:flow} with
reference to~\eqref{eq:ode_system}, which we now particularize
to~\eqref{eq:ode_system_ce}. By Lemma~\ref{lem:HS}, $\flow$ is of
class $\C2$.  Define
\begin{equation}
  \label{eq:18}
  (q^\flat, p^\flat) (t)
  \coloneqq
  \flow (t, 0, \sqrt2)
  \qquad \mbox{ and } \qquad
  (q^\sharp, p^\sharp) (t)
  \coloneqq
  \flow (t, 0, 2) \,,
\end{equation}
$q^\flat$, respectively $q^\sharp$, being the leftmost, respectively
rightmost, red line in Figure~\ref{fig:fc2}.

\begin{lemma}
  \label{lmm:ce5}
  Let $H$ be as in~\eqref{eq:11} and $u_o$ be as
  in~\eqref{eq:12}. Define the set
  \begin{equation}
    \label{eq:17}
    D \coloneqq
    \left([0, +\infty\mathclose[ \times \{2\}\right)
    \cup
    \left(\{0\} \times \mathopen]0,2]\right) .
  \end{equation}
  Then, there exists a unique map
  \begin{equation}
    \label{eq:19}
    \fonction{\Delta}{\mathopen]0, +\infty\mathclose[ \times \mathopen]0, +\infty\mathclose[}{D}{(t, x)}{(q_o, p_o)}
  \end{equation}
  such that
  \begin{equation}
    \label{eq:20}
    \flow_q(t, q_o, p_o) = x
    \quad \text{and} \quad
    \forall \, s \in \mathopen]0, t\mathclose[, \; \flow_q(s, q_o, p_o) > 0 \,.
  \end{equation}
  Moreover,
  \begin{enumerate}[\rm (1)]
  \item\label{item:20} $\Delta$ is continuous.
  \item\label{item:21} $\Delta$ is monotone, in the sense that setting
    $\Delta (t_o,x_o) = (0,p_o)$ and $\Delta (t_o,x_o') = (0,p_o')$,
    if $0 < x_o < x_o' < q^\sharp(t_o)$, then $p_o < p_o'$.
  \item\label{item:22} For all
    $x \in \mathopen]0, +\infty\mathclose[$,
    $\lim_{t\to 0+} \Delta (t,x) = (x,2)$.
  \end{enumerate}

\end{lemma}

\begin{proofof}{Lemma~\ref{lmm:ce5}}
  We split the proof in several steps.
  \paragraph{For all $(t,x) \in \mathopen]0, +\infty\mathclose[^2$,
    there exists $(p_o,q_o) \in D$ satisfying~\eqref{eq:20}.} Fix
  $(t,x) \in \mathopen]0, +\infty\mathclose[ \times \mathopen]0,
  +\infty\mathclose[$. If $x = q^\sharp(t)$ as in~\eqref{eq:18}, then
  set $(q_o, p_o) = (0, 2)$. Otherwise, introduce the functions
  \begin{displaymath}
    \fonction{h}{[0, +\infty\mathclose[}{\R}{q_o}{\flow_q(t, q_o, 2) - x}
    \quad \mbox{ and } \quad
    \fonction{k}{\mathopen]0, 2]}{\R}{p_o}{\flow_q(t, 0, p_o) - x\,.}
  \end{displaymath}
  Note that if $x > q^\sharp(t)$ then $h (0) < 0$ and
  $h (x+1) = \flow_q (t,x+1,2) - x > 0$ by Lemma~\ref{lmm:ce1}. By the
  Intermediate Value Theorem, there exists a $q_o$ such that
  $h (q_o) = 0$, hence $\flow_q (t,q_o,2) = x$. By
  Lemma~\ref{lmm:ce1}, for all $s > 0$, $\flow_q (s,q_o,2) > 0$,
  proving~\eqref{eq:20} in the case $x > q^\sharp(t)$.

  If $x \in [q^\flat (t), q^\sharp (t)]$, then
  $k (\sqrt2) \leq 0 \leq k (2)$ by~\eqref{eq:18}. By the Intermediate
  Value Theorem, there exists $p_o$ such that $k (p_o) = 0$, i.e.,
  $\flow_q (t,0,p_o) = x$. Then, by Lemma~\ref{cor:ce4}, the right
  part of~\eqref{eq:20} follows in the case
  $x \in [q^\flat (t), q^\sharp (t)]$.

  Similarly, if $x \in \mathopen]0, q^\flat(t)\mathclose[$, then
  $k (\sqrt2) >0$ and $\lim_{\check p \to 0+} k (\check p) = -x <
  0$. By the Intermediate Value Theorem, we can define
  \begin{displaymath}
    p_o
    \coloneqq
    \max \left\{ \pi_o \in [0,\sqrt2] \colon k (\pi_o) = 0\right\} \,.
  \end{displaymath}
  Hence, for all $\check p \in \mathopen] p_o, \sqrt2\mathclose[$,
  $k (\check p) >0$. Proceed now by contradiction: assume there exists
  $s \in \mathopen]0,t\mathclose[$ such that $\flow_q (s,0,p_o ) <
  0$. By Lemma~\ref{lmm:reallyLast}, we get $t > \mathcal{T}
  (p_o)$. Using Lemma~\ref{lem:Last} and the Intermediate Value
  Theorem, it follows that there exists
  $p_o' \in \mathopen]p_o, \sqrt2\mathclose[$ such that
  $\mathcal{T} (p_o') = t$. By Lemma~\ref{lmm:reallyLast}, this
  implies that $\flow_q (t,0, p_o') = 0$ and therefore $k (p_o') < 0$,
  which contradicts the choice of $p_o$.

  \paragraph{$\Delta$ is uniquely defined.} For all
  $(t,x) \in \mathopen]0, +\infty\mathclose[ \times \mathopen]0,
  +\infty\mathclose[$ the uniqueness of a $(q_o,p_o)$
  satisfying~\eqref{eq:20} follows from the monotonicity properties
  proved above. Indeed, recalling $q^\sharp$ as defined
  in~\eqref{eq:18}, if $0 < q_o < \tilde q_o$ and
  $p_o = \tilde p_o = 2$, then, by Lemma~\ref{lmm:ce2}, for all
  $s \in [0, +\infty\mathclose[$,
  $ q^\sharp (s) < \flow_q(s, q_o, p_o) < \flow_q(s, \tilde q_o,
  \tilde p_o)$. On the other hand, if $q_o = \tilde q_o =0$ and
  $0 < p_o < \tilde p_o < 2$, then by Lemma~\ref{lmm:ce3},
  Lemma~\ref{lmm:reallyLast} and Lemma~\ref{cor:ce4}, for all $s$ such
  that $\flow_q (\tau,q_o,p_o) \geq 0$ for all $\tau \in [0,s]$, we
  have
  $\flow_q(s, q_o, p_o) < \flow_q(s, \tilde q_o, \tilde p_o) <
  q^\sharp (s)$. Finally, if $q_o = 0$, $\tilde q_o >0$,
  $p_o \in \mathopen]0,2\mathclose[$ and $\tilde p_o =2$, then
  Lemma~\ref{cor:ce4} and Lemma~\ref{lmm:ce2} ensure that for all $s$
  such that $\flow_q (\tau,q_o,p_o) \geq 0$ for all $\tau \in [0,s]$,
  we have
  $\flow_q(s, q_o, p_o) < q^\sharp (s) < \flow_q(s, \tilde q_o, \tilde
  p_o)$. The uniqueness of $(q_o,p_o)$ follows.

  \paragraph{$\Delta$ is continuous.} For any $(q_o,p_o) \in D$ and
  $(t,x) \in \R_+^2$, if $(q_o,p_o) = \Delta (t,x)$ then
  by~\eqref{eq:conservation_ce}, we have
  $\modulo{\flow_p (t,q_o,p_o)} \leq \sqrt{p_o^2 + 2}$, so that
  by~\eqref{eq:ode_system_ce},
  $\modulo{\flow_q (t,q_o,p_o) - q_o} \leq t \sqrt{p_o^2 +
    2}$. Therefore,
  \begin{equation}
    \label{eq:29}
    \modulo{q_o} \leq x + t \sqrt{p_o^2 + 2} \,.
  \end{equation}
  Choose now a sequence $(t_n,x_n)$ in
  $\mathopen]0, +\infty \mathclose[ \times \mathopen]0, +\infty
  \mathclose[$ converging to $(t,x)$ also in
  $\mathopen]0, +\infty \mathclose[ \times \mathopen]0, +\infty
  \mathclose[$. Define $(q_o^n, p_o^n) = \Delta (t_n,x_n)$. The
  sequence $p_o^n$ is in $[0,2]$ by~\eqref{eq:17}
  and~\eqref{eq:19}. By~\eqref{eq:29}, also the sequence $q_o^n$ is
  bounded, since also $(t_n,x_n)$ is bounded. Call $(q_o,p_o)$ the
  limit of any convergent subsequence, so that
  $(q_o,p_o) \in \overline{D}$. By the continuity of $\flow$ proved in
  Lemma~\ref{lem:HS}. up to a subsequence we have
  \begin{equation}
    \label{eq:22}
    \flow_q (t, q_o, p_o)
    =
    \lim_{n\to+\infty} \flow_q (t_n, q_o^n, p_o^n)
    =
    \lim_{n\to+\infty} x_n
    = x \,.
  \end{equation}
  This also shows that $(q_o,p_o) \in D$. Otherwise, if
  $(q_o,p_o) \in \overline{D} \setminus D$, then $(q_o,p_o) = (0,0)$
  and for all $t \in \R$, $\flow_q (t,0,0) = 0 \neq x$.

  Since $(q_o^n, p_o^n) = \Delta (t_n,x_n)$, then
  $x_n = \flow_q (t_n,q_o^n,p_o^n)$. Thus, if
  $s \in \mathopen]0, t_n\mathclose[$, then
  $\flow_q (s,q_o^n,p_o^n) > 0$. In the limit $n \to +\infty$, we have
  $x = \flow_q (t,q_o,p_o)$ and if $s \in \mathopen]0, t\mathclose[$,
  then $\flow_q (s,q_o,p_o) \geq 0$.

  The possible behaviors of $s \to \flow_q(s, q_o, p_o)$ classified in
  Lemma~\ref{lmm:ce1}, Lemma~\ref{lmm:ce3} and in
  Lemma~\ref{lmm:reallyLast} ensure that for all
  $s \in \mathopen]0, t\mathclose[$ we have $\flow_q(s, q_o, p_o) > 0$
  so that also the second condition in~\eqref{eq:20} is met and
  $\Delta (t,x) = (q_o,p_o)$, the limit $(q_o,p_o)$ being independent
  of the subsequence. This completes the proof of the continuity of
  $\Delta$.

  \paragraph{Proof of~\eqref{item:21} and~\eqref{item:22}.} Fix a
  positive $x$. Let $t_n$ be any positive sequence converging to
  $0$. Then, $\flow_q\left(t_n, \Delta (t_n,x) \right) = x$. The
  bound~\eqref{eq:29} ensures that, up to a subsequence,
  $\lim_{n \to +\infty} \Delta (t_n,x) =\xi$, with
  $\xi \in \overline{D}$ satisfying $\flow_q (0, \xi) = x$. Hence,
  $\xi = (x,2)$, proving~\eqref{item:22}.

  The monotonicity of $\Delta$ follows from Lemma~\ref{cor:ce4},
  completing the proof of~\eqref{item:21}.
\end{proofof}

\begin{proposition}
  \label{prop:ce6}
  Let $H$ be as in~\eqref{eq:11} and $u_o$ be as
  in~\eqref{eq:12}. Recall the notations~\eqref{eq:flow}
  and~\eqref{eq:19}. The function
  \begin{equation}
    \label{eq:21}
    \fonction{u}{\mathopen]0, +\infty\mathclose[ \times (\R \setminus\{0\})}{\R}{(t,x)}{\left\{
        \begin{array}{ccl}
          \flow_p\left(t, \Delta(t, x)\right) & \text{if} & x > 0 \,,
          \\
          -\flow_p\left(t, \Delta(t, -x)\right) & \text{if} & x < 0 \,.
        \end{array}
      \right.}
  \end{equation}
  is in $\L{\infty}(\mathopen]0,+\infty\mathclose[ \times \R;\R)$,
  solves~\eqref{eq:cl} with datum~\eqref{eq:12} in the sense of
  Definition~\ref{def:entropy_solution}, it is a classical strong
  solution outside $x=0$ and outside $|x| = q^\sharp (t)$, it is
  continuous along $|x| = q^\sharp (t)$ and there is an entropic
  stationary shock along $x=0$ for $t > \pi / (2\sqrt2)$.
\end{proposition}

The lack of differentiability along $|x| = q^\sharp (t)$ is visible in
Figure~\ref{fig:nice}.

\begin{proofof}{Proposition~\ref{prop:ce6}}
  Call $\Gamma$ the graph of the map $t \mapsto q^\sharp (t)$ as
  defined in~\eqref{eq:18} and define
  $\Omega \coloneqq \mathopen]0, +\infty\mathclose[^2 \setminus
  \Gamma$. Note that by Lemma~\ref{lmm:ce5} and~\eqref{eq:21},
  $u \in \C0(\mathopen]0, +\infty\mathclose[ \times (\R
  \setminus\{0\}); \R)$.

  \paragraph{Claim~1: $u \in \C{1}(\Omega;\R)$ and is a classical
    solution to~\eqref{eq:cl}--\eqref{eq:11}--\eqref{eq:12} in
    $\Omega$.} This follows from an application of the Implicit
  Function Theorem. Indeed, let $(t_o, x_o) \in \Omega$. Then, either
  $\Delta(t_o, x_o) = (q_o, 2)$ or $\Delta (t_o,x_o) = (0, p_o)$ for
  suitable $q_o > 0$ or $p_o \in \mathopen]0,2\mathclose[$. Thus,
  $\flow_q(t_o, q_o, 2) - x_o = 0$ or $\flow_q (t_o,0,p_o) - x_o = 0$
  .  Introduce the functions
  \begin{displaymath}
    \fonction{F}{\mathopen]0, +\infty\mathclose[^3}{\R}{(t, x, q)}{\flow_q(t, q, 2) - x}
    \quad
    \fonction{G}{\mathopen]0, +\infty\mathclose[^2 \times \mathopen]0, 2\mathclose[}{\R}{(t, x, p)}{\flow_q(t, 0, p) - x \,.}
  \end{displaymath}
  By Lemma~\ref{lem:HS}, both $F$ and $G$ are of class $\C{2}$,
  $F(t_o, x_o, q_o) = 0$ or $G (t_o, x_o, p_o) = 0$.

  Moreover, Lemma~\ref{lmm:ce2} implies that for all $t > 0$,
  $q \mapsto \flow_q(t, q, 2)$ is increasing and therefore,
  \begin{displaymath}
    \forall \, (t, q) \in \mathopen]0,+\infty\mathclose[^2,
    \quad \partial_q \flow_q (t, q, 2) \geq 0 \,.
  \end{displaymath}
  If $\partial_q \flow_q (t_o, q_o, 2) = 0$, then $t_o$ minimizes the
  map $t \mapsto \partial_q \flow_q (t, q_o, 2)$ so that
  $\frac{\d{~}}{\d{t}} \partial_q \flow_q (t, q_o, 2) = 0$ and
  $y \colon t \mapsto \partial_q \flow_q (t, q_o, 2)$ solves the
  Cauchy problem
  \begin{displaymath}
    \left\{
      \begin{aligned}
        \ddot y(t) & = -g''\left(\flow_q(t, q_o, 2)\right) y(t) \\
        y(t_o) & = 0
        \\
        \dot y(t_o) & = 0 \,.
      \end{aligned}
    \right.
  \end{displaymath}
  The uniqueness of solutions is ensured by Cauchy Lipschitz Theorem,
  we thus have that $y \equiv 0$. On the other hand,
  deriving~\eqref{eq:ode_system_ce} with respect to $q_o$, we see that
  $y$ also solves
  \begin{displaymath}
    \left\{
      \begin{aligned}
        \ddot y(t) & = -g''\left(\flow_q(t, q_o, 2)\right) y(t) \\
        y(0) & = 1
        \\
        \dot y(0) & = 0 \,,
      \end{aligned}
    \right.
  \end{displaymath}
  which is a contradiction. Therefore,
  $\partial_q F(t_o, x_o, q_o) > 0$.

  Similarly, Lemma~\ref{cor:ce4} implies that for all $t>0$,
  $p \mapsto \flow_q (t,x,p)$ is increasing and therefore
  \begin{displaymath}
    \forall \, (t, p) \in \mathopen]0,+\infty\mathclose[ \times \mathopen]0, 2\mathclose[ \,,
    \quad \partial_p \flow_q (t, 0, p) \geq 0 \,.
  \end{displaymath}
  If $\partial_p \flow_q (t_o, 0, p_o) = 0$, then $t_o$ minimizes the
  map $t \mapsto \partial_p \flow_q (t, 0, p_o)$ so that
  $\frac{\d{~}}{\d{t}} \partial_p \flow_q (t, 0, p_o) = 0$ and
  $t \mapsto \partial_p \flow_q (t, 0, p_o)$ solves the Cauchy problem
  \begin{displaymath}
    \left\{
      \begin{aligned}
        \ddot y(t) & = -g''\left(\flow_q(t, 0, p_o)\right) y(t) \\
        y(t_o) & = 0
        \\
        \dot y(t_o) & = 0\,.
      \end{aligned}
    \right.
  \end{displaymath}
  The uniqueness of solutions ensured by Cauchy Lipschitz Theorem, we
  thus have that $y \equiv 0$.  On the other hand,
  deriving~\eqref{eq:ode_system_ce} with respect to $p_o$, we see that
  $y$ also solves
  \begin{displaymath}
    \left\{
      \begin{aligned}
        \ddot y(t) & = -g''\left(\flow_q(t, 0, p_o)\right) y(t) \\
        y(0) & = 0
        \\
        \dot y(0) & = 1 \,,
      \end{aligned}
    \right.
  \end{displaymath}
  which is a contradiction. Therefore,
  $\partial_p G(t_o, x_o, p_o) > 0$.

  The Implicit Function Theorem allows us to obtain a locally unique
  map $Q$ such that $q_o = Q(t_o, x_o)$ from the relation
  $F(t_o, q_o, x_o) = 0$ and, in the same way, to obtain
  $p_o = P(t_o,x_o)$ from the relation $G (t_o,x_o,p_o) = 0$, with
  both functions $Q$ and $P$ of class $\C1$. Note that
  by~\eqref{eq:19}, by~\eqref{eq:20}, by~\eqref{item:20} in
  Lemma~\ref{lmm:ce5} and by the local uniqueness of $Q$ and $P$, we
  get
  \begin{displaymath}
    \Delta (t,x) =
    \left\{
      \begin{array}{lr@{\,}c@{\,}l}
        \left(Q (t,x),2\right)
        & \mbox{if }x
        & >
        & q^\sharp (t)
        \\
        \left(0, P (t,x)\right)
        & \mbox{if }x
        & <
        & q^\sharp (t)
      \end{array}
    \right.
  \end{displaymath}
  and, by~\eqref{eq:21}, the $\C1$ regularity of $u$ in $\Omega$ is
  proved.

  We now prove that $u$ solves~\eqref{eq:cl} with $H$ as
  in~\eqref{eq:11} and initial datum~\eqref{eq:12}. To this aim,
  observe that the map $x \mapsto u (t,x)$ is odd, for all
  $t \in \R_+$.

  Assume $x>0$. Then, by the Implicit Function Theorem and
  by~\eqref{eq:21}, for all $(t, x) \in \Omega$ we have
  \begin{displaymath}
    u (t,x)
    =
    \left\{
      \begin{array}{lr@{\,}c@{\,}l}
        \flow_p\left(t, Q (t,x),2\right)
        & \mbox{if }x
        & >
        & q^\sharp (t)
        \\
        \flow_p\left(t, 0, P (t,x)\right)
        & \mbox{if }x
        & <
        & q^\sharp (t)
      \end{array}
    \right.
    \mbox{ and } \;
    \begin{array}{r@{\,}c@{\,}l@{\qquad}r@{\,}c@{\,}l}
      \partial_t Q
      & =
      & - \frac{\partial_t \flow_q}{\partial_q \flow_q}
        = -\frac{\flow_p}{\partial_q \flow_q}
      & \partial_x Q
      & =
      & \frac{1}{\partial_q \flow_q} \,;
      \\
      \partial_t P
      & =
      & -\frac{\partial_t \flow_q}{\partial_p \flow_q}
        = -\frac{\flow_p}{\partial_p \flow_q}
      & \partial_x P
      & =
      & \frac{1}{\partial_p \flow_q} \,.
    \end{array}
  \end{displaymath}
  Hence, recalling also~\eqref{eq:ode_system_ce}
  \begin{displaymath}
    \partial_t u + u \; \partial_x u
    =
    \left\{
      \begin{array}{l}
        \partial_t \flow_p
        +
        \partial_q \flow_p \; \partial_t Q
        + \flow_p \; \partial_q \flow_p \; \partial_x Q
        =
        - g'
        \\
        \partial_t \flow_p
        + \partial_p \flow_p \; \partial_t P
        + \flow_p \; \partial_p \flow_p \; \partial_x P
        =
        - g'
      \end{array}
    \right.
  \end{displaymath}
  ensuring that $u$ solves~\eqref{eq:cl}--\eqref{eq:11} in the
  classical sense in $\Omega$. The case $x<0$ is entirely similar
  by~\eqref{eq:21}, since $H$ is even in $x$ and $p$.

  \paragraph{Claim~2: Conclusion.}  The monotonicity proved in
  Lemma~\ref{lmm:ce5} ensures that $\Delta$, and hence $u$, admits
  traces along $x=0$ for all $t>0$ and $u (t,0-) = - u (t,0+)$,
  because $u (t)$ is odd. Since $p \mapsto \frac{p^2}{2}$ is an even
  function, by~\eqref{eq:11} and~\eqref{eq:21}
  $H\left(0, u(t, 0+)\right) = H\left(0, u(t, 0-)\right)$. Hence,
  either $u (t, 0+) = u (t, 0-)$, or Rankine--Hugoniot conditions hold
  along the stationary discontinuity along $x=0$.

  Assume that $u (t, 0+) \neq u(t, 0-)$. Then, $u (t, 0+) \neq
  0$. Moreover, by~\eqref{eq:21}, for a positive sequence $x_n$
  converging to $0$, we have that $u (t, x_n) \neq 0$ and has a fixed
  signed for all $n$. Hence, $\Delta (t,x_n) = (0, p_n)$ with
  $p_n \neq 0$. Lemma~\ref{lmm:ce3} and Lemma~\ref{lmm:reallyLast}
  then ensure that $p_n \in \mathopen]0, \sqrt2 \mathclose[$. Up to a
  subsequence, $\lim_{n\to +\infty}p_n = p_*$ for a suitable $p_* >0$
  (for, otherwise, $u$ would vanish).
  \begin{flalign*}
    u (t,0+) %
    & = \lim_{n\to+\infty} u (t, x_n)
    \\
    & = \lim_{n\to+\infty} \flow_p \left(t, \Delta (t,x_n)\right) %
    & [\mbox{By~\eqref{eq:21}}]
    \\
    & = \lim_{n\to+\infty} \flow_p (t, 0, p_n) %
    & [\mbox{By~\eqref{eq:19}--\eqref{eq:20}}] %
    \\
    & = \lim_{n\to+\infty} \frac{\d{~}}{\d{t}} \flow_q (t,0,p_n) %
    & [\mbox{By~\eqref{eq:ode_system_ce}}]
    \\
    & = \frac{\d{~}}{\d{t}} \flow_q (t,0,p_*) \,. %
    & [\mbox{By Lemma~\ref{lem:HS}}]
  \end{flalign*}
  Note that $\flow_q (t,0,p_*) = 0$, so that by~\eqref{eq:20} the map
  $t \mapsto \flow_q (t,0,p_*)$ passes from positive to negative at
  $t$, showing that $\frac{\d{~}}{\d{t}} \flow_q (t,0,p_*) \leq 0$, so
  that $u (t,0+) \leq u (t,0-)$. By~\eqref{eq:convexity}, we obtain
  the Lax Entropy Inequality~\cite[Section~11.9]{DafermosBook} at
  $x=0$.

  Consider now the initial datum. By~\eqref{item:22} in
  Lemma~\ref{lmm:ce5}, for any
  $x \in \mathopen]0, +\infty\mathclose[$,
  $\lim_{t \to 0+} u (t,x) = \lim_{t\to 0+} \flow_p\left(t, \Delta
    (t,x)\right) = \flow_p\left(0, x, 2\right) = 2 = u_o (x)$. The
  case $x<0$ is entirely analogous.

  Along $\Gamma$, $u$ is continuous so that Rankine--Hugoniot and Lax
  entropy conditions are met. Hence, $u$ is an entropy solution
  to~\eqref{eq:cl}--\eqref{eq:11}--\eqref{eq:12} both where $x>0$ and,
  by symmetry, also where $x<0$. Along $x=0$, Rankine--Hugoniot
  conditions and the usual Lax entropy inequalities are met, both if
  $u$ is continuous or not. As $t\to 0+$, $u (t)$ pointwise converges
  to the initial datum~\eqref{eq:12}. A standard argument then ensures
  that $u$ solves~\eqref{eq:cl}--\eqref{eq:11}--\eqref{eq:12} in the
  sense of Definition~\ref{def:entropy_solution}, see~\cite[Conclusion
  in the proof of Lemma~3.3]{CPS2022} for the details. The proof is
  completed.
\end{proofof}

\begin{proofof}{Theorem~\ref{th:counter_example}}
  Let $T > 0$. Set $w = u(T, \cdot)$ with $u$ as in~\eqref{eq:21}. By
  Proposition~\ref{prop:ce6}, $I_T^{CL}(w) \neq \emptyset$ and
  $u_o \in I_T^{CL}(w)$, where $u_o$ is as in~\eqref{eq:12}. Moreover,
  thanks to the regularity of $u$ outside $x=0$ proved in
  Proposition~\ref{prop:hypOk} and
  to~\cite[Theorem~4.1]{Dafermos1977}, since $\pi_{w}(\R) = \R$, by
  Theorem~\ref{th:belonging} we deduce that $I_T^{CL}(w) = \{u_o\}$.

  To complete the proof, use $\mathcal{T}$ as defined
  in~\eqref{eq:15}, note that a shock in $u$ first arises at time
  $\inf_{p_o \in [0, \sqrt2]} \mathcal{T} (p_o) / 2 = \pi / (2\sqrt
  2)$ by~\eqref{item:26} in Lemma~\ref{lem:Last}. The growth of the
  shock size follows from the fact that $\mathcal{T}$ is strictly
  increasing.
\end{proofof}

\appendix

\section{Appendix}
\label{sec:appendix}

\begin{lemma}
  \label{lem:Sturm}
  The polynomial
  $P (x) \coloneqq x^8 - \frac{32}{7} x^6 + \frac{59}{7} x^4 - 8x^2 -
  \frac{6}{7}$ is negative for all $x \in [-1,1]$.
\end{lemma}

\begin{proofof}{Lemma~\ref{lem:Sturm}}
  The Sturm sequence, see~\cite{Sturm1829}, of $P$ is:
  \begin{displaymath}
    \begin{array}{r@{\qquad}c@{\qquad}c}
      & \mbox{sign at }-1
      & \mbox{sign at } 1
      \\[6pt]
      \displaystyle
      x^8 - \frac{32}{7} x^6 + \frac{59}{7} x^4 - 8x^2 - \frac{6}{7}
      & -
      & -
      \\[6pt]
      \displaystyle
      8x^7 - \frac{192}{7} x^5 + \frac{236}{7} x^3 - 16x
      & +
      & -
      \\[6pt]
      \displaystyle
      \frac{8}{7} x^6 - \frac{59}{14} x^4 + 6x^2 + \frac{6}{7}
      & +
      & +
      \\[6pt]
      \displaystyle
      -\frac{29}{14} x^5 + \frac{58}{7} x^3 + 22x
      & -
      & +
      \\[6pt]
      \displaystyle
      -\frac{5}{14} x^4 - \frac{526}{29} x^2 - \frac{6}{7}
      & -
      & -
      \\[6pt]
      \displaystyle
      -\frac{3972}{35} x^3 - \frac{944}{35} x
      & +
      & -
      \\[6pt]
      \displaystyle
      \frac{3639116}{201579} x^2 + \frac{6}{7}
      & +
      & +
      \\[6pt]
      \displaystyle
      \frac{137451770}{6368453} x
      & -
      & +
      \\[6pt]
      \displaystyle
      -\frac{6}{7}
      & -
      & -
    \end{array}
  \end{displaymath}
  Sturm Theorem, see~\cite{Sturm1829}, ensures that $P$ has
  $4 - 4 = 0$ roots in $[-1, 1]$. Therefore, for all $x \in [-1,1]$,
  $P (x)$ has the sign of $P(0) = -6/7 < 0$.
\end{proofof}

\begin{lemma}
  \label{lem:NoHope}
  Let $\mathcal{A}$ be as in~\eqref{eq:26} with $g$ as
  in~\eqref{eq:11}. Then, $\mathcal{A}$ is continuous on
  $\mathopen]0,1\mathclose[$.
\end{lemma}

\begin{proofof}{Lemma~\ref{lem:NoHope}}
  For
  $(r,\theta) \in \mathopen]0,1 \mathclose[ \times [ 0,1\mathclose[$,
  define
  $\mathcal{B} (r,\theta) \coloneqq \left. r \middle/ \sqrt{g (r) - g
      (\theta\, r)}\right.$. $\mathcal{B}$ is positive. For any
  $\epsilon \in \mathopen]0, 1/2\mathclose[$, fix
  $r \in [\epsilon, 1-\epsilon]$. Then, for $\theta \in [0,1/2]$,
  $\mathcal{B} (r,\theta) \leq \max_{[\epsilon, 1-\epsilon] \times
    [0,1/2]} \mathcal{B}$, while for $\theta \in [1/2, 1\mathclose[$
  \begin{displaymath}
    \mathcal{B} (r,\theta)
    \leq
    \dfrac{1}{\sqrt{1-\theta} \sqrt{\min_{\rho \in [\epsilon/2, 1-\epsilon]} g' (\rho)}} \,.
  \end{displaymath}
  Hence, $\mathcal{B}$ is continuous and dominated, therefore
  $\mathcal{A}$ is continuous, too.
\end{proofof}

\smallskip

\paragraph{Acknowledgment:}
The first author was partly supported by the GNAMPA~2022 project
\emph{Evolution Equations, Well Posedness, Control and Applications}.
This research was funded, in whole or in part, by l’Agence Nationale
de la Recherche (ANR), project ANR-22-CE40-0010. For the purpose of
open access, the second author has applied a CC-BY public copyright
license to any Author Accepted Manuscript (AAM) version arising from
this submission.

{\small

  \bibliography{hID}

\begin{thebibliography}{10}

\bibitem{AdimurthiGhoshal2020}
{Adimurthi} and S.~S. Ghoshal.
\newblock Exact and optimal controllability for scalar conservation laws with
  discontinuous flux, 2020.
\newblock Preprint.

\bibitem{MR4188826}
F.~Ancona and M.~T. Chiri.
\newblock Attainable profiles for conservation laws with flux function
  spatially discontinuous at a single point.
\newblock {\em ESAIM Control Optim. Calc. Var.}, 26:Paper No. 124, 33, 2020.

\bibitem{zbMATH06067762}
F.~Ancona, O.~Glass, and K.~T. Nguyen.
\newblock Lower compactness estimates for scalar balance laws.
\newblock {\em Commun. Pure Appl. Math.}, 65(9):1303--1329, 2012.

\bibitem{MR1616586}
F.~Ancona and A.~Marson.
\newblock On the attainable set for scalar nonlinear conservation laws with
  boundary control.
\newblock {\em SIAM J. Control Optim.}, 36(1):290--312, 1998.

\bibitem{MR3394696}
B.~Andreianov, C.~Donadello, S.~S. Ghoshal, and U.~Razafison.
\newblock On the attainable set for a class of triangular systems of
  conservation laws.
\newblock {\em J. Evol. Equ.}, 15(3):503--532, 2015.

\bibitem{BCDBook}
M.~Bardi and I.~Capuzzo-Dolcetta.
\newblock {\em Optimal control and viscosity solutions of
  {H}amilton--{J}acobi-Bellman equations}.
\newblock Springer Science \& Business Media, 2008.

\bibitem{BarlesBook}
G.~Barles.
\newblock An introduction to the theory of viscosity solutions for first-order
  {H}amilton--{J}acobi equations and applications.
\newblock In {\em Hamilton-Jacobi equations: approximations, numerical analysis
  and applications}, pages 49--109. Springer, 2013.

\bibitem{MR1722801}
E.~N. Barron, P.~Cannarsa, R.~Jensen, and C.~Sinestrari.
\newblock Regularity of {H}amilton-{J}acobi equations when forward is backward.
\newblock {\em Indiana Univ. Math. J.}, 48(2):385--409, 1999.

\bibitem{MR2347697}
A.~Bressan and B.~Piccoli.
\newblock {\em Introduction to the mathematical theory of control}, volume~2 of
  {\em AIMS Series on Applied Mathematics}.
\newblock AIMS, Springfield, MO, 2007.

\bibitem{CSBook}
P.~Cannarsa and C.~Sinestrari.
\newblock {\em Semiconcave functions, {H}amilton--{J}acobi equations, and
  optimal control}, volume~58.
\newblock Springer Science \& Business Media, 2004.

\bibitem{Chicone1987}
C.~Chicone.
\newblock The monotonicity of the period function for planar {H}amiltonian
  vector fields.
\newblock {\em Journal of Differential equations}, \textbf{69}(3):310--321,
  1987.

\bibitem{MR2795714}
S.~Cifani and E.~R. Jakobsen.
\newblock Entropy solution theory for fractional degenerate
  convection-diffusion equations.
\newblock {\em Ann. Inst. H. Poincar\'{e} C Anal. Non Lin\'{e}aire},
  28(3):413--441, 2011.

\bibitem{ClarkeBook}
F.~Clarke.
\newblock {\em {Functional Analysis, Calculus of Variations and Optimal
  Control}}.
\newblock Graduate Texts in Mathematics, Vol. 264. {Springer}, 2013.

\bibitem{MR2374224}
G.~M. Coclite and N.~H. Risebro.
\newblock Viscosity solutions of {H}amilton-{J}acobi equations with
  discontinuous coefficients.
\newblock {\em J. Hyperbolic Differ. Equ.}, 4(4):771--795, 2007.

\bibitem{CP2020}
R.~M. Colombo and V.~Perrollaz.
\newblock Initial data identification in conservation laws and
  {H}amilton--{J}acobi equations.
\newblock {\em J. Math. Pures et Appl.}, \textbf{138}:{1--27}, 2020.

\bibitem{CPS2022}
R.~M. Colombo, V.~Perrollaz, and A.~Sylla.
\newblock Conservation laws and {H}amilton--{J}acobi equations with space
  inhomogeneity, 2022.
\newblock Preprint.

\bibitem{procHYP2022}
R.~M. Colombo, V.~Perrollaz, and A.~Sylla.
\newblock Peculiarities of space dependent conservation laws: Inverse design
  and asymptotics, 2023.
\newblock Preprint.

\bibitem{MR3489384}
M.~Corghi and A.~Marson.
\newblock On the attainable set for scalar balance laws with distributed
  control.
\newblock {\em ESAIM Control Optim. Calc. Var.}, 22(1):236--266, 2016.

\bibitem{CL1983}
M.~G. Crandall and P.-L. Lions.
\newblock Viscosity solutions of {H}amilton--{J}acobi equations.
\newblock {\em Transactions of the American Mathematical Society},
  \textbf{277}(1):1--42, 1983.

\bibitem{Dafermos1977}
C.~M. Dafermos.
\newblock Generalized characteristics and the structure of solutions of
  hyperbolic conservation laws.
\newblock {\em Indiana University Mathematics Journal},
  \textbf{26}(6):1097--1119, 1977.

\bibitem{DafermosBook}
C.~M. Dafermos.
\newblock {\em Hyperbolic conservation laws in continuum physics}, volume 325
  of {\em Grundlehren der Mathematischen Wissenschaften}.
\newblock Springer-Verlag, Berlin, fourth edition, 2016.

\bibitem{zbMATH02183964}
C.~de~Lellis and F.~Golse.
\newblock A quantitative compactness estimate for scalar conservation laws.
\newblock {\em Commun. Pure Appl. Math.}, 58(7):989--998, 2005.

\bibitem{EZ2020}
C.~Esteve and E.~Zuazua.
\newblock The inverse problem for hamilton--jacobi equations and semiconcave
  envelopes.
\newblock {\em SIAM Journal on Mathematical Analysis},
  \textbf{52}(6):5627--5657, 2020.

\bibitem{MR4503821}
C.~Esteve-Yag\"{u}e and E.~Zuazua.
\newblock Reachable set for {H}amilton-{J}acobi equations with non-smooth
  {H}amiltonian and scalar conservation laws.
\newblock {\em Nonlinear Anal.}, 227:Paper No. 113167, 2023.

\bibitem{EvansBook}
L.~C. Evans.
\newblock {\em Partial Differential Equations}.
\newblock American Mathematical Society, Providence, R.I., 2010.

\bibitem{MR47234}
E.~Hopf.
\newblock The partial differential equation {$u_t+uu_x=\mu u_{xx}$}.
\newblock {\em Comm. Pure Appl. Math.}, 3:201--230, 1950.

\bibitem{KarlsenRisebro2002}
K.~H. Karlsen and N.~H. Risebro.
\newblock A note on front tracking and equivalence between viscosity solutions
  of {H}amilton-{J}acobi equations and entropy solutions of scalar conservation
  laws.
\newblock {\em Nonlinear Anal.}, 50(4, Ser. A: Theory Methods):455--469, 2002.

\bibitem{Kruzhkov1970}
S.~N. Kruzhkov.
\newblock First order quasilinear equations with several independent variables.
\newblock {\em Mathematics of the USSR-Sbornik}, \textbf{81}(123):228--255,
  1970.

\bibitem{Lax}
P.~D. Lax.
\newblock Hyperbolic systems of conservation laws. {I}{I}.
\newblock {\em Comm. Pure Appl. Math.}, 10:537--566, 1957.

\bibitem{LZ2021}
T.~Liard and E.~Zuazua.
\newblock Initial data identification for the one-dimensional {B}urgers
  equation.
\newblock {\em IEEE Transactions on Automatic Control}, 2021.

\bibitem{Sturm1829}
J.~C.~F. Sturm.
\newblock Mémoire sur la résolution des équations numériques.
\newblock {\em Bulletin des Sciences de Férussac}, \textbf{11}:419--425, 1829.

\bibitem{SyllaPhD}
A.~Sylla.
\newblock {\em {Heterogeneity in scalar conservation laws: approximation and
  applications}}.
\newblock Thesis, {Universit{\'e} de Tours}, July 2021.
\newblock \url{https://hal.archives-ouvertes.fr/tel-03303049}.

\end{thebibliography}

  \bibliographystyle{abbrv}

}

\end{document}